\documentclass[10.0pt]{amsart}
\usepackage{graphicx} 
\usepackage{tikz}
\usepackage{epsfig}
\usepackage{amsmath, amssymb, euscript,  enumerate}
\usepackage{pxfonts}
\usepackage{bbm}
\usepackage{color}
\usepackage{verbatim}

\usepackage{scalerel}[2014/03/10]
\usepackage[usestackEOL]{stackengine}
\def\intavg{\,\ThisStyle{\ensurestackMath{%
      \stackinset{c}{0\LMpt}{c}{0\LMpt}{\SavedStyle-}
{\SavedStyle\phantom{\int}}}%
       \setbox0=\hbox{$\SavedStyle\int\,$}\kern-\wd0}\int}

\newcommand{\supp}{\text{\rm supp}}

\newcommand{\ap}{\alpha}             
\newcommand{\bt}{\beta}
\newcommand{\gm}{\gamma}             \newcommand{\Gm}{\Gamma}
\newcommand{\dt}{\delta}             
\newcommand{\vep}{\varepsilon}
\newcommand{\zt}{\zeta}
\newcommand{\Th}{\Theta}

\newcommand{\ld}{\lambda}            \newcommand{\Ld}{\Lambda}
\newcommand{\sm}{\sigma}             
\newcommand{\vp}{\varphi}
             \newcommand{\Om}{\Omega}
\newcommand{\vr}{\varrho}            \newcommand{\iy}{\infty}
            
\newcommand{\f}{\frac}

\newcommand{\fA}{{\mathfrak A}}
\newcommand{\fB}{{\mathfrak B}}
\newcommand{\fC}{{\mathfrak C}}
\newcommand{\fD}{{\mathfrak D}}
\newcommand{\fE}{{\mathfrak E}}

\newcommand{\fa}{{\mathfrak a}}

\newcommand{\BN}{{\mathbb N}}

\newcommand{\BR}{{\mathbb R}}

\newcommand{\cA}{{\mathcal A}}
\newcommand{\cB}{{\mathcal B}}
\newcommand{\cC}{{\mathcal C}}

\newcommand{\cF}{{\mathcal F}}

\newcommand{\cI}{{\mathcal I}}
\newcommand{\cJ}{{\mathcal J}}

\newcommand{\cL}{{\mathcal L}}

\newcommand{\cT}{{\mathcal T}}

\newcommand{\la}{\langle}          \newcommand{\ra}{\rangle}
\newcommand{\s}{\setminus}         \newcommand{\ep}{\epsilon}
\newcommand{\n}{\nabla}            \newcommand{\e}{\eta}

    \newcommand{\ds}{\displaystyle}
\newcommand{\ls}{\lesssim } \newcommand{\gs}{\gtrsim}

 \newcommand{\pf }{\noindent{\it Proof. }}

\newcommand{\aee }{\text{\rm a.e.}} \newcommand{\diam }{\text{\rm diam}}
  \newcommand{\pv }{\text{\rm p.v.}}

\newcommand{\rX}{{\text{\rm X}}}

\newcommand{\loc}{{\text{\rm loc}}}

\newcommand{\bmo }{{\text{BMO}}}
\newcommand{\ei}{\text{\rm ess inf }}

\newtheorem{thm}[subsection]{Theorem}
\newtheorem{lemma}[subsection]{Lemma}
\newtheorem{cor}[subsection]{Corollary}

\newtheorem{prop}[subsection]{Proposition}

\newtheorem{defn}[subsection]{Definition}

\numberwithin{equation}{section}

\title[Nonlocal Harnack inequalities for Nonlocal double phase equations I ]{Nonlocal Harnack inequalities for \\ Nonlocal double phase equations I ;\\
with nonnegative bounded modulating  coefficient \\ with no H\"older condition }
\author{ Yong-Cheol Kim  }
\begin{document}
\begin{abstract}
In this paper, by applying the De Giorgi-Nash-Moser theory we prove nonlocal Harnack inequalities for (locally nonnegative in $\Om$) weak solutions to nolocal double phase equations
\begin{equation*}\begin{cases}\cL u =0 & \text{ in $\Om$,} \\
u=g & \text{ in $\BR^n\s\Om$ }
\end{cases}\end{equation*}
where $\Om\subset\BR^n$ ($n\ge 2$) is a bounded domain with Lipschitz boundary, $\cL$ is the nonlocal double phase operator $\cL$ given by 
\begin{equation*}\begin{split}\cL u(x)=&\pv\int_{\BR^n}|u(x)-u(y)|^{p-2}(u(x)-u(y))K_{ps}(x,y)\,dy \\
&+\pv\int_{\BR^n}\fa(x,y)|u(x)-u(y)|^{q-2}(u(x)-u(y))K_{qt}(x,y)\,dy,
\end{split}
\end{equation*} $\fa(x,y)=\fa(y,x)\ge 0$, $\|\fa\|_{L^\iy(\BR^n\times\BR^n)}<\iy$, $\ei\fa>0$, and 
$ps\ge qt$ for $0<s,t<1$ and $1<p\le q<\iy$. 

In addition, we get local boundedness with explicit formula and weak Harnack inequalities for their weak supersolutions.
\end{abstract}

\address{$\bullet$ Yong-Cheol Kim : Department of Mathematics Education, Korea University, Seoul 02841, Republic of Korea   $\,\&\,$ School of Mathematics, Korea Institute for Advanced Study, Seoul 02455, Republic of Korea  }

\email{ychkim@korea.ac.kr}

\maketitle

\tableofcontents

\section{ Introduction }
Double phase equations appear in natural environment with the energy density switching between two different types of degenerate components, according to the characteristic of a {\it modulating coefficient} $\fa$ which determines the {\it phase}. The energy functional of the double phase equation has the nonstandard growth conditions. In the local theory, their energy functionals are given by
\begin{equation}\fE(u)=\int_{\Om}\bigl(|Du|^p+\fa(x)|Du|^q\bigr)\,dx
\end{equation}
where $1<p\le q$, $\Om\subset\BR^n\,(n\ge 2)$ is a bounded domain with Lipschitz boundary and $\fa:\Om\to\BR$ is a nonnegative bounded function. when $p=q$, we call them functionals with {\it standard growth conditions}, and also when $p<q$, we call them functionals with {\it nonstandard growth conditions}. Initially, the low order regularity theory for the minimizers of the functionals (1.1) 
was formulated and studied by Marcellini \cite{M1,M2}. Nowadays, the regularity results of the minimizers are well-understood in the local realm of energy functionals of double phase equations with nonstandard growth \cite{G, Ma, Mi}. Moreover, these functionals can be generalized as
\begin{equation*}\fE(u)=\int_\Om F(x,u,Du)\,dx
\end{equation*} heir
where $F:\Om\times\BR\times\BR^n\to[0,\iy)$ is a Carath\'eodory function satisfying 
\begin{equation*}|z|^p\ls F(x,t,z)\ls|z|^q+1,\,\,1<p<q,\,(x,t,z)\in\Om\times\BR\times\BR^n.
\end{equation*}

We are now interested in studying a nonlocal version of double phase equations mentioned in the above.
For our aim, we consider the following nonlocal double phase operator
\begin{equation*}\begin{split}\cL u(x)=&\pv\int_{\BR^n}|u(x)-u(y)|^{p-2}(u(x)-u(y))K_{ps}(x,y)\,dy \\
&+\pv\int_{\BR^n}\fa(x,y)|u(x)-u(y)|^{q-2}(u(x)-u(y))K_{qt}(x,y)\,dy
\end{split}
\end{equation*}
where $K_{h}:\BR^n\times\BR^n\to\BR$ are symmetric kernels satisfying 
\begin{equation}\begin{split}
\f{\ld}{|x-y|^{n+h}}&\le K_{h}(x,y)=K_{h}(y,x)\le\f{\Ld}{|x-y|^{n+h}},\,\,0<\ld\le\Ld<\iy,\,h>0\\
\end{split}
\end{equation} for some $s,t,p,q$ with
\begin{equation}qt\le ps,\,\,0<s,t<1<p\le q<\iy\end{equation}
and $\fa:\BR^n\times\BR^n\to\BR$ is the modulating coefficient satisfying
\begin{equation}\fa(x,y)=\fa(y,x)\ge 0,\,\, \|\fa\|_{L^\iy(\BR^n\times\BR^n)}:=\|\fa\|_\iy<\iy\,\,\text{ and }\,\, \ei\fa>0.
\end{equation}
For $\tau\in\BR$ and $p>1$, we set $H_p(\tau)=|\tau|^{p-2}\tau$. Then the nonlocal double phase operator $\cL$ in the above can be written as
\begin{equation*}\cL u(x)=\cL_p u(x)+\cL_q^\fa u(x)
\end{equation*} where the operators $\cL_p$ and $\cL_q^\fa$ are given by
\begin{equation*}\begin{split}
\cL_p u(x)&=\pv\int_{\BR^n}H_p(u(x)-u(y))K_{ps}(x,y)\,dy, \\
\cL_q^\fa u(x)&=\pv\int_{\BR^n}\fa(x,y)H_q(u(x)-u(y))K_{qt}(x,y)\,dy.
\end{split}\end{equation*}
For $x_0\in\Om$ and $r>0$ with $B_r(x_0)\subset\Om$, we define the nonlocal tail $\cT_r(u;x_0)$ of the function $u$ related with $\cL\,$ by
\begin{equation*}\cT_r(u;x_0)=r^{ps}\int_{\BR^n\s B_r(x_0)}\biggl(\,\f{|u(y)|^{p-1}}{|y-x_0|^{n+ps}}+\|\fa\|_{\iy}\f{|u(y)|^{q-1}}{|y-x_0|^{n+qt}}\biggr)\,dy.
\end{equation*}

\,

\noindent{\bf Notations.} We introduce the notations briefly for the readers as follows.

\,\noindent$\bullet$ For $r>0$ and $x_0\in\BR^n$, let us denote by $B_r^0=B_r(x_0)$, $B_r=B_r(0)$.

\,\noindent$\bullet$ For two quantities $a$ and $b$, we write $a\lesssim b$ (resp.
$a\gtrsim b$) if there is a universal constant $C>0$ ({\it depending
only on $\ld,\Ld,n,p,q,s,\Om$ and $\fa$}$)$  such that $a\le C\,b$ (resp. $b\le
C\,a$).

\,\noindent$\bullet$ For $a,b\in\BR$, we denote by 
$a\vee b=\max\{a,b\}$ and $a\wedge b=\min\{a,b\}.$

\,\noindent$\bullet$ Let $\cF(\BR^n)$ be the set of all real-valued
Lebesgue measurable functions on $\BR^n$.

\,\noindent$\bullet$ For $u\in C(B^0_r)$, we consider the norm 
$\|u\|_{C(B^0_r)}=\sup_{x\in B^0_r}|u(x)|.$ For $\gm\in(0,1)$, the {\it $\gm^{th}$ H{$\ddot o$}lder seminorm} of $u$ on $B^0_r$ is defined by
$$[u]_{C^{\gm}(B^0_r)}=\sup_{x,y\in B^0_r,\,x\neq y}\f{|u(x)-u(y)|}{|x-y|^{\gm}}
$$ and the {\it $\gm^{th}$ H{$\ddot o$}lder norm} of $u$ on $B^0_r$ is defined by 
$\|u\|_{C^{\gm}(B^0_r)}=\|u\|_{C(B^0_r)}+[u]_{C^{\gm}(B^0_r)}.$

\,\noindent$\bullet$ The {\it essential infimum} of the above-mentioned modulating coefficient $\fa$ is defined by
$$\ei\fa=\sup_{c\in\BR}\bigl\{c:|\{(x,y)\in\BR^n\times\BR^n:\fa(x,y)<c\}|_{2n}=0\bigr\},$$
where $|\cdot|_{2n}$ is the Lebesgue measure on $\BR^n\times\BR^n$. Also, we denote by
$$\|\fa\|_*:=\f{\ei\fa}{\|\fa\|_\iy}.$$
Then we see that $0<\|\fa\|_*\le 1$, if $\ei\fa>0$; here we note that $\|\fa\|_*=1$ if $\fa$ is nonzero constant almost everywhere.

\,

In this paper, our aim is to study the Dirichlet  problems related with the above-mentioned nonlocal double phase operator $\cL$ as follows;
\begin{equation}\begin{cases}\cL u =0 & \text{ in $\Om$,} \\
u=g & \text{ in $\BR^n\s\Om$ }
\end{cases}\end{equation}
where $\Om\subset\BR^n$ is a bounded domain with Lipchitz boundary. In particular, we are concerned with nonlocal Harnack inequalities for (locally nonnegative in $\Om$)  weak solutions of the nonlocal double phase equations (1.5).

The research about local double phase equations with nonnegative bounded modulating coefficients has been ongoing by numerous researchers. By using the frozen functional method, Baroni-Colombo-Mingione \cite{BCM} obtained Harnack inequalities for weak solutions to local double phase equations with nonnegative bounded H\"older continuous modulating coefficients.

In the literature, nonlocal double phase equations with nonnegative bounded modulating coefficients $\fa$ were first studied by De Filippis and Palatucci \cite{DP} who obtained H\"older regularity for bounded viscosity solutions to the nonlocal double phase equations. Local H\"older continuity for bounded weak solutions of the nonlocal double phase equations with the condition $qt\le ps$ was obtained in \cite{FZ}. Under H\"older condition of $\fa$,  local H\"older continuity for weak solutions to the nonlocal double phase equations with the condition $ps\le qt$ and their local boundedness with no explicit formula was shown in \cite{BOS}. Also, in case of $ps\ge qt$, H\"older regularity for weak solutions to nonlocal double phase equations equipped only with nonnegative bounded modulating coefficient with strictly positive essential infimum and with no H\"older condition was obtained in \cite{K2}.

For nonlocal Harnack estimates, it is one of the main steps to get the local boundedness with explicit formula. We realize the local boundedness with explicit upper bound in Theorem 1.3 below. When $ps\ge qt$, we obtain nonlocal Harnack inequalities for weak solutions to the nonlocal double phase equations equipped with nonnegative bounded modulating coefficient $\fa$ with strictly positive essential infimum and  with no H\"older condition. 
To our knowledge, this result is the first one for nonlocal Harnack inequality of their weak solutions under the assumption.
But we need a restriction $\ei\fa>0$ in $\BR^n\times\BR^n$. Notwithstanding the restriction, the infimum and the supremum of $\fa$ on $\BR^n\times\BR^n$ could vanish and become infinite, respectively. Along this line, we will study the case $ps\le qt$ in the future.

In what follows, via fractional Sobolev inequality and Morrey inequality, for $p>1$ we define the {\it fractional Sobolev conjugate} $p_s^*$ of $p$ by
\begin{equation*}p_s^*=\begin{cases} \ds\f{np}{n-ps} &\text{ if $\,ps<n$,} \\
\text{ arbitrary number in $[p,\iy)$ } &\text{ if $\,ps\ge n.$ }\end{cases}
\end{equation*}
We note that this implies that, for any $p>1$,
\begin{equation*}p\le q\le p_s^*\,\,\,\,\Leftrightarrow\,\,\,\,\begin{cases}\,\, p\le q\le\ds\f{np}{n-ps} &\text{ if $\,ps<n$, }\\
\,\,p\le q<\iy &\text{ if $\,ps\ge n$. }
\end{cases}
\end{equation*}

\,

We now state our main result which is called {\it nonlocal Harnack inequalities} for (locally nonnegative in $\Om$) weak solutions to the nonlocal double phase equations (1.5) as follows. 

\begin{thm} Suppose that $qt\le ps\le n$ for $s,t\in(0,1)$ and $1<p\le q\le p_s^*$, and let $g\in L^{q-1}_{qt}(\BR^n)$. 
If $u\in\fD_g(\Om)$ is a weak solution of the nonlocal double phase equation $(1.5)$ with $u\ge 0$ in $B^0_R\subset\Om$, then there is a constant $C>0$ depending only on $n,s,p,q,\ld,\Ld,\Om$ and $\fa$ such that
\begin{equation*}\begin{split}
\sup_{B^0_r}u\le C\,\biggl[\inf_{B^0_r}u+\biggl(\f{r}{R}\biggr)^{\f{ps}{p-1}}
\cT^{\f{1}{p-1}}_R(u^-;x_0)\biggr]\,\,\,\,\text{ for any $r\in(0,R)$. } \\
\end{split}\end{equation*}
\end{thm}

Theorem 1.1 implies the classical Harnack inequalities for globally nonnegative weak solutions to the nonlocal double phase equations (1.5) as follows.

\begin{cor} Suppose that $qt\le ps\le n$ for $s,t\in(0,1)$ and $1<p\le q\le p_s^*$, and let $g\in L^{q-1}_{qt}(\BR^n)$. 
If $u\in\fD_g(\Om)$ is a weak solution of the nonlocal double phase equation $(1.5)$ with $u\ge 0$ in $\BR^n$, then there is a constant $C>0$ depending only on $n,s,p,q,\ld,\Ld,\Om$ and $\fa$ such that
\begin{equation*}\begin{split}
\sup_{B^0_r}u\le C\,\inf_{B^0_r}u
\end{split}\end{equation*}
\end{cor}

We next state one of our results which is called the {\it local boundedness} of weak subsolutions of the nonlocal double phase equations (1.5) as follows. 

\begin{thm} Let $qt\le ps\le n$ for $s,t\in(0,1)$ and $1<p\le q\le p_s^*$, $g\in L^{q-1}_{qt}(\BR^n)$ and $B^0_{2r}\subset\Om$.  If $u\in\fD_g^-(\Om)$ is a weak subsolution to the nonlocal double phase equation $(1.5)$, then there is a constant $C_0>0$ depending only on $n,s,p,q,\ld,\Ld,\diam(\Om)$ and $\fa$ such that
\begin{equation}\begin{split}
\sup_{B^0_r}u&\le\dt\,\cT^{\f{1}{p-1}}_r(u^+;x_0)
+C_0\,h(\dt)\,\biggl(\,\intavg_{B^0_{2r}}\Theta(u^+)\,dx\biggr)^{\f{1}{p}}
\end{split}\end{equation} 
for any $\dt\in(0,1]$, 
where the function $\Theta$ is given by $\Theta(z):=z^p+\|a\|_{\iy}\,z^q,\,\,z\ge 0,$ and 
$$h(\dt)=\begin{cases} \,\,\dt^{-\f{(p-1)n}{sp^2}} & \text{ for $\,\,ps<n$, }\\
\,\,\dt^{-\f{q(p-1)}{p(q-p)}} & \text{ for $\,\,ps=n$.}
\end{cases}$$
\end{thm} 

The following {\it logarithmic estimate} plays an important role in proving H\"older continuity of weak solutions and the infimum estimate of weak supersolutions. This estimate makes it possible to show that the logarithm of weak supersolutions to the nonlocal double phase equation (1.5) is a function with locally bounded mean oscillation.

\begin{thm} Suppose that $qt\le ps\,$ for $s,t\in(0,1)$ and $1<p\le q<\iy$ and $g\in L^{q-1}_{qt}(\BR^n)$. 
If $u\in\fD^+_g(\Om)$ is a weak supersolution of the nonlocal double phase equation $(1.5)$ with $u\ge 0$ in $B^0_R\subset\Om$, then we have that
$$\iint_{B^0_r\times B^0_r}\biggl|\ln\biggl(\f{u(x)+b}{u(y)+b}\biggr)\biggr|^p\,\f{dx\,dy}{|x-y|^{n+pt_{p}}}
\ls r^{n-pt_{p}}\biggl[1+\f{1}{b^{p-1}}\biggl(\f{r}{R}\biggr)^{ps}\,\cT_R(u^-;x_0)\biggr].$$
for any $b>0$ and any $r\in(0,R/2)$, where $t_{p}:=t(1-1/p)\in(0,1)$.
\end{thm}

It is well-known in \cite{T} that, if $v\in W^{s,p}(\Om)$ for $ps>n$ with $0<s<1<p<\iy$ where $B\subset\Om$ is a ball, then there is a universal constant $c=c(n,s,p)>0$ such that 
\begin{equation*}[v]_{C^{s-\f{n}{p}}(B)}\le c\,[v]_{W^{s,p}(B)},
\end{equation*}
which is called the {\it fractional Morrey inequality}. Combining Theorem 1.4 with this inequality, we can easily derive the following corollary which furnishes a sort of nonlocal Harnack inequality for its weak supersolution and classical Harnack inequality for its weak solution.

\begin{cor} Suppose that $qt\le ps\,$ for $s,t\in(0,1)$ and $1<p\le q<\iy$ and $g\in L^{q-1}_{qt}(\BR^n)$. 
If $\,p>1+n/t\,$ and $u\in\fD^+_g(\Om)$ is a weak supersolution of the nonlocal double phase equation $(1.5)$  with $u\ge 0$ in $B^0_R\subset\Om$, then for any $b>0$, any $r\in(0,R/2)$ and any $x,y\in B^0_r$ we have that
\begin{equation*}\exp\biggl(-\biggl[1+\f{1}{b^{p-1}}\biggl(\f{r}{R}\biggr)^{ps}\,\cT_R(u^-;x_0)\biggr]\biggr)
\ls\f{u(x)}{u(y)}\ls\exp\biggl[1+\f{1}{b^{p-1}}\biggl(\f{r}{R}\biggr)^{ps}\,\cT_R(u^-;x_0)\biggr].
\end{equation*}
Moreover, if $\,p>1+n/t\,$ and $u\in\fD_g(\Om)$ is a weak supersolution of the nonlocal double phase equation $(1.5)$ with $u\ge 0$ in $\BR^n$, then  for any $r\in(0,R/2)$ and any $x,y\in B^0_r$ we have that
\begin{equation*}1/e\ls\f{u(x)}{u(y)}\ls e.
\end{equation*}
\end{cor}

For the proof of Theorem 1.1, we need the following {\it nonlocal weak Harnack inequalities} whose proof will be executed in Section 8 from a version of the Krylov-Safonov covering theorem and the infimum estimate of weak supersolution of the nonlocal double phase equation (1.5).

\begin{thm} Let $qt\le ps\le n$ for $0<s,t<1$ and $1<p\le q\le p_s^*$, $g\in L^{q-1}_{qt}(\BR^n)$. 
Suppose that $u\in\fD^+_g(\Om)$ is a weak supersolution of the nonlocal double phase equation $(1.5)$ with $u\ge 0$ in $B^0_R\subset\Om$. Then there exists some $h_0\in(0,1)$ such that 
\begin{equation*}\sup_{h\in(0,h_0]}\biggl(\,\intavg_{B^0_{2r}}u^h\,dx\biggr)^{\f{1}{h}}\le 2^{\f{1}{h_0}}\biggl[\inf_{B^0_r}u+\f{1}{2}\biggl(\f{r}{R}\biggr)^{\f{ps}{p-1}}\,\cT^{\f{1}{p-1}}_R(u^-;x_0)\biggr]
\end{equation*} for any $r\in(0,R/2).$
\end{thm}


\begin{figure}
    \centering
    \begin{tikzpicture}[scale=1.2, every node/.style={font=\tiny}]
    \clip (-1.5,-0.6) rectangle (5, 5.5);

    \draw[->] (-0.7,0) -- (4,0) node[right] {$p$};
    \draw[->] (0,-0.3) -- (0,4.5) node[above] {$q$};

    \fill[gray!30] (1,1) -- (1,1.25) -- (2.3, {1.25*2.3}) -- (2.3, 2.3) -- cycle;

    \fill[gray!30] (3,3) -- (3, {1.25*3}) -- (3.5, {1.25*3.5}) -- (3.5, 3.5) -- cycle;

    \fill[red!30] (2.3,2.3) -- (2.3, {1.25*2.3}) -- (3, {1.25*3}) -- (3, 3) -- cycle;

    \draw[domain=0:3.5, red, thick] plot (\x, \x);
    \node[red] at (3.5,3) {$q=p$};

    \draw[domain=0:3.5, blue, thick] plot (\x, {1.25*\x});
    \node[blue] at (3,4.5) {$q=\frac{s}{t}p$};

    \draw[domain=0:1.5, green!60!black, thick, samples=100] plot (\x, {-5.29/(\x-2.3)-2.3});
    \node[green!60!black, right] at (0.4, 4.5) {$q=\frac{np}{n-ps}$};

    \draw[dotted] (0, 23/13) -- (1, 23/13);
    \node[left=1pt] at (0, 23/13) {$\frac{n}{n-s}$};

    \draw[dotted] (0,1.25) -- (1,1.25);
    \node[left=1pt] at (0,1.25) {$\frac{s}{t}$};

    \draw[dotted] (1,0) -- (1, 23/13);
    \node[below=1pt] at (1,0) {$1$};

    \draw[-] (2.3,0) -- (2.3, 4.5); 
    \node[below=1pt] at (2.3,0) {$\frac{n}{s}$};

    \draw[dotted] (3,0) -- (3, 3.75);
    \node[below=1pt] at (3,0) {$1+\frac{n}{t}$};

    \end{tikzpicture}
    \caption{$ps<n$, $\,\frac{s}{t}\le\frac{n}{n-s}$}
\end{figure}

\begin{figure}
    \centering
    \begin{tikzpicture}[scale=0.9, xscale=0.85, yscale=0.4, every node/.style={font=\tiny}]
    
    \clip (-2,-2) rectangle (7.5, 16);

    \fill[gray!30, samples=100] (1,1) -- (1, 4/3) 
        plot[domain=1:2.4] (\x, {-16/(\x-4)-4}) 
        -- (4, 10)  
        -- (4, 4)  
        -- (1, 1)   
        -- cycle;

    \fill[red!30] (4,4) -- (4, 10) -- (5, 12.5) -- (5, 5) -- cycle;

    \fill[gray!30] (5,5) -- (5, 12.5) -- (6, 6*2.5) -- (6, 6) -- cycle;

    \draw[->] (-1,0) -- (7,0) node[right] {$p$};
    \draw[->] (0,-1) -- (0,15) node[above] {$q$};

    \draw[domain=0:6, red, thick] plot (\x, \x);
    \node[red] at (5.8, 5) {$q=p$};

    \draw[domain=0:6, blue, thick] plot (\x, {2.5*\x});
    \node[blue, left] at (5.8, 15) {$q=\frac{s}{t}p$};

    \draw[domain=0:3, green!60!black, thick, samples=100] plot (\x, {-16/(\x-4)-4});
    \node[green!60!black, right, font=\tiny] at (1.5, 12.7) {$q=\frac{np}{n-ps}$};

    \draw[dotted] (2.4, 0) -- (2.4, 6);
    \node[below=1pt] at (2.4, 0) {$\frac{n(s-t)}{s^2}$};

    \draw[dotted] (0, 6) -- (2.4, 6);
    \node[left=1pt] at (0, 6) {$\frac{n(s-t)}{st}$};
    
    \draw[dotted] (0,4/3) -- (1,4/3);
    \node[left=1pt] at (0,4/3) {$\frac{n}{n-s}$};
    
    \draw[dotted] (0,2.5) -- (1,2.5);
    \node[left=1pt] at (0,2.5) {$\frac{s}{t}$};

    \draw[dotted] (1,0) -- (1,2.5);
    \node[below=1pt] at (1,0) {$1$};

    \draw[-] (4,0) -- (4, 15);
    \node[below=1pt] at (4,0) {$\frac{n}{s}$};

    \draw[dotted] (5,0) -- (5, 12.5);
    \node[below=1pt] at (5,0) {$1+\frac{n}{t}$};

    \end{tikzpicture}
    \caption{$ps<n$, $\,\frac{s}{t}>\frac{n}{n-s}$}
\end{figure}

\begin{figure}
    \centering
    \begin{tikzpicture}[scale=1.2, every node/.style={font=\tiny}]
    \clip (-1.5,-0.6) rectangle (5, 5.5);

    \draw[->] (-0.7,0) -- (4,0) node[right] {$p$};
    \draw[->] (0,-0.3) -- (0,4.5) node[above] {$q$};

    \fill[gray!30] (1,1) -- (1,1.25) -- (2.3, {1.25*2.3}) -- (2.3, 2.3) -- cycle;

    \fill[gray!30] (3,3) -- (3, {1.25*3}) -- (3.5, {1.25*3.5}) -- (3.5, 3.5) -- cycle;

    \fill[red!30] (2.3,2.3) -- (2.3, {1.25*2.3}) -- (3, {1.25*3}) -- (3, 3) -- cycle;

    \draw[domain=0:3.5, red, thick] plot (\x, \x);
    \node[red] at (3.5,3) {$q=p$};

    \draw[domain=0:3.5, blue, thick] plot (\x, {1.25*\x});
    \node[blue] at (3,4.5) {$q=\frac{s}{t}p$};




    \draw[dotted] (1,0) -- (1, 23/13);
    \node[below=1pt] at (1,0) {$1$};

    \draw[-] (2.3,0) -- (2.3, 4.5); 
    \node[below=1pt] at (2.3,0) {$\frac{n}{s}$};

    \draw[dotted] (3,0) -- (3, 3.75);
    \node[below=1pt] at (3,0) {$1+\frac{n}{t}$};

    \end{tikzpicture}
    \caption{$ps=n$}
\end{figure}


{\bf Remark.} From Theorem 1.1 and Corollary 1.5, we obtain nonlocal Harnack inequalities for weak solutions to the nonlocal double phase equation (1.5) as you can see in the gray zone of Figure 1, Figure 2 and Figure 3. For some more details, in Figure 1 and Figure 2 we give the gray region in which nonlocal Harnack inequalities was obtained in case of $ps<n$ with two possible cases, and also in Figure 3 we furnish it in case of $ps=n$. However, we could not fill up the red zone of Figure 1, Figure 2 and Figure 3. So it will be very interesting to ask whether the nonlocal Harnack inequalities work or not in the red zone.

\,

The article is organized as follows. In Section 2, we give the function spaces and the definition of weak solutions of the nonlocal double phase equation given in (1.5). In Section 3, we furnish several elemetary lemmas and well-known lemmas which are very useful in applying the De Giorgi-Nash-Moser method. In Section 4, we give embedding properties of the fractional Sobolev spaces.
In Section 5, we obtain a sort of {\it nonlocal Caccioppoli type inequality} and several useful local properties of weak solutions to the nonlocal double phase equation. In Section 6, we prove that the logarithm of weak solutions to the nonlocal double phase equation is a function with locally bounded mean oscillation. In Section 7, we shall obtain the supremum and infimum estimate of weak supersolutions to the nonlocal double phase equation. In Section 8, we prove nonlocal weak Harnack inequalities for weak solutions to the nonlocal double phase equation by applying the Krylov-Safonov covering theorem and the results obtained in Section 7. In Section 9, we finally obtain nonlocal Harnack inequalities by applying the local boundedness and nonlocal weak Harnack inequalities of its weak solutions.

\section{ Preliminaries }

Let $\Om\subset\BR^n$ be a bounded domain with Lipschitz boundary and let $0<s<1$. For $p>1$, let $\rX^{s,p}(\Om)$ be the linear function space of all Lebesgue measurable functions $u\in\cF(\BR^n)$ on $\BR^n$ such that $u|_\Om\in L^p(\Om)$ and
\begin{equation}\iint_{\BR^{2n}_\Om}\f{|u(x)-u(y)|^p}{|x-y|^{n+ps}}\,dx\,dy<\iy\end{equation}
where $\BR^{2n}_S:=(\BR^n\times\BR^n)\s[(\BR^n\s S)\times(\BR^n\s S)]$ for a set $S\subset\BR^n$. Then we see that $\rX^{s,p}(\Om)$ is a normed space with the norm $\|\cdot\|_{\rX^{s,p}(\Om)}$ given by 
\begin{equation*}\|u\|_{\rX^{s,p}(\Om)}=\|u\|_{L^p(\Om)}+\biggl(\,\iint_{\BR^{2n}_\Om}\f{|u(x)-u(y)|^p}{|x-y|^{n+ps}}\,dx\,dy\biggr)^{\f{1}{p}}<\iy\,\,\text{ for $u\in\rX^{s,p}(\Om)$.}
\end{equation*}
For $s\in(0,1)$, $p\in[1,\iy)$ and an open set $U\subset\BR^n$, the {\it fractional Sobolev space} $W^{s,p}(U)$
is defined as the set of all functions $u\in L^p(U)$ satisfying that
\begin{equation}\begin{split}
\|u\|_{W^{s,p}(U)}&:=\|u\|_{L^p(U)}+[u]_{W^{s,p}(U)} \\
&:=\|u\|_{L^p(U)}+\biggl(\,\iint_{U\times U}\f{|u(x)-u(y)|^p}{|x-y|^{n+ps}}\,dx\,dy\biggr)^{\f{1}{p}}<\iy,
\end{split}\end{equation}
that is, an intermediary Banach space between $L^p(U)$ and $W^{1,p}(U)$ with the norm $\|\cdot\|_{W^{s,p}(U)}$.
Also, we denote by $W^{s,p}_0(U)$ the closure of $C^\iy_c(U)$ in the norm $\|\cdot\|_{W^{s,p}(U)}$.
When $U=\BR^n$ in (2.2), similarly we define the space $W^{s,p}(\BR^n)$ for $p\ge 1$ and $s\in(0,1)$.

If $ps<n$ for $p\ge 1$ and $s\in(0,1)$, then it is well-known \cite{DPV} that there is a constant $c=c(n,p,s,\Om)>0$ such that
\begin{equation}\|f\|_{L^\tau(\Om)}\le c\,\|f\|_{W^{s,p}(\Om)}
\end{equation} for any $f\in W^{s,p}(\Om)$ and $\tau\in[p,p_s^*]$ where $p_s^*=\f{pn}{n-ps}$.
Moreover, there is a constant $c=c(n,p,s)>0$ such that
\begin{equation}\|f\|_{L^{p_s^*}(\BR^n)}\le c\,[f]_{W^{s,p}(\BR^n)}\,\,\text{ for any $f\in W^{s,p}(\BR^n)$ } 
\end{equation}
and 
\begin{equation}\|f\|_{L^\tau(\BR^n)}\le c\,\|f\|_{W^{s,p}(\BR^n)}\,\,\text{ for any $f\in W^{s,p}(\BR^n)$ and $\tau\in[p,p_s^*]$. }
\end{equation}
For $1<p\le q<\iy$ and $0<s\le t<1$, let $\Th_*$ be the function defined by
\begin{equation*}\Th_*(x,y;z):=\f{z^p}{|x-y|^{ps}}+\fa(x,y)\f{z^q}{|x-y|^{qt}}, \,\,\,x,y\in\BR^n,\,\,z\ge 0.
\end{equation*}
Then we consider two the double integral values given by
\begin{equation*}\fB(v;E):=\iint_{\BR^{2n}_E}\Th_*(x,y;|v(x)-v(y)|)\,\f{dx\,dy}{|x-y|^n},
\end{equation*}
where 
$v\in\cF(\BR^n)$ and $\BR^{2n}_E$ is a set defined in (2.1).
Now we define a function space $\fD(\Om)$ related with weak solutions to the Dirichlet problem (1.5) by
\begin{equation*}\fD(\Om)=\{v\in\cF(\BR^n):v|_\Om\in L^p(\Om)\,\text{ and }\,\fB(v;\Om)<\iy\}.
\end{equation*}
Then it is obvious that $\fB(v;E)<\iy$ for $v\in\fD(\Om)$ and 
$\fD(\Om)\subset\rX^{s,p}(\Om)\subset W^{s,p}(\Om).$
Also we set 
\begin{equation*}\fD_0(\Om):=\{v\in\fD(\Om):u=0\,\text{ $\aee$ in $\BR^n\s\Om$ }\}.
\end{equation*}

We consider the tail space for the nonlocal double phase equation mentioned in (1.5), which reflects the nonlocal effect of the nonlocal equation. The tail space $L^{q-1}_{qt}(\BR^n)$ is defined by
\begin{equation*}L^{q-1}_{qt}(\BR^n)=\bigl\{v\in L^{q-1}_{\loc}(\BR^n):\int_{\BR^n}\f{|v(y)|^{q-1}}{(1+|y|)^{n+qt}}\,dy<\iy\bigr\}.
\end{equation*}
When $ps\ge qt$, it is easy to check that 
\begin{equation}W^{s,p}(\BR^n)\subset L^{p-1}_{ps}(\BR^n),\,\,\,\,W^{t,q}(\BR^n)\subset L^{q-1}_{qt}(\BR^n)
\subset L^{p-1}_{ps}(\BR^n)\,\,
\end{equation} and
\begin{equation}v\in L^{q-1}_{qt}(\BR^n)\, \text{ if and only if $\,\cT_r(v;x_0)<\iy\,$ for any $r\in(0,\iy)$. }
\end{equation}
Moreover, we easily see from the fractional Sobolev inequality that 
\begin{equation*}\begin{split}\fD(\Om)\subset L^\tau(\Om)\,\text{ for }\,
\begin{cases}\,p\le \tau\le\ds\f{np}{n-ps} &\text{ when $ps<n$, } \\
\,p\le \tau<\iy&\text{ when $ps\ge n$. } \end{cases}
\end{split}\end{equation*}
For $g\in L^{q-1}_{qt}(\BR^n)$, we consider the convex subset of $\fD(\Om)$ given by
\begin{equation*}\fD_g(\Om)=\{v\in\fD(\Om):g-v\in\fD_0(\Om)\}.
\end{equation*}
Note that $\fD_g(\Om)=\fD_0(\Om)$ whenever $g=0$. In order to define supersolutions and subsolutions to the nonlocal double phase equation (1.5), we define 
\begin{equation*}\fD^\pm(\Om)=\{v\in\fD(\Om):(g-v)_\pm\in\fD_0(\Om)\}
\end{equation*} and we consider a bilinear form $\la\cdot,\cdot\ra_\cL:\fD(\Om)\times\fD(\Om)\to\BR$ defined by
\begin{equation*}\begin{split}\la u,\vp\ra_\cL=
\iint_{\BR^n\times\BR^n}&\biggl(H_p(u(x)-u(y))(\vp(x)-\vp(y))K_{ps}(x,y)  \\
&\qquad+\fa(x,y)H_q(u(x)-u(y))(\vp(x)-\vp(y))K_{qt}(x,y)\biggr)\,dx\,dy.
\end{split}
\end{equation*}

\begin{defn} Let $qt\le ps\,$ for $s,t\in(0,1)$ and $1<p\le q<\iy$ and $g\in L^{q-1}_{qt}(\BR^n)$. Then we say that a function $u\in\fD_g^-(\Om)$ $($ $\fD_g^+(\Om)$ $)$ is a {\rm weak subsolution (weak supersolution)} to the nonlocal double phase equation $(1.5)$, if it satisfies that
\begin{equation*}\begin{split}
\la u,\vp\ra_\cL\le (\, \ge \,) \,\,0
\end{split}\end{equation*} for every nonnegative $\vp\in\fD_0(\Om)$. Also, we say that a function $u\in\fD_g(\Om)$ is a {\rm weak solution} to the nonlocal double phase equation $(1.5),\,$ if it is both a weak subsolution and a weak supersolution to the equation.
\end{defn}

\section{ Several useful elementary inequalities }
In this section, we mention and prove several elementary inequalities which are useful in solving our problems.

\begin{lemma} Let $\vartheta_1,\vartheta_2:(0,\iy)\to(0,\iy)$ be functions satisfying that $\vartheta_1\le\vartheta_2$ and
\begin{equation}\lim_{z\to 0^+}\f{\vartheta_i(z)}{z}\ge 0\,\,\text{ and }\,\,\f{\vartheta_i(z)}{z}\text{ is increasing on $(0,\iy)\,$ for $i=1,2$. }
\end{equation}
If $a,b\in(0,\iy)$, then we have the following inequality
\begin{equation}\f{b}{a}\le \f{\vartheta_2(b)}{\vartheta_1(a)}+\mathbbm{1}_{\{a\ge b\}}(a,b)
\end{equation} where $\mathbbm{1}_{\{a\ge b\}}:(0,\iy)\times(0,\iy)\to\{0,1\}$ is the function given by $\mathbbm{1}_{\{a\ge b\}}(a,b)=1$ if $a\ge b$ and $\mathbbm{1}_{\{a\ge b\}}(a,b)=0$ if $a<b$.
\end{lemma}

\pf Assume that $a\ge b$. Then (3.2) is equivalent to
$$1\le\f{\vartheta_2(b)/b}{\vartheta_1(a)/a}+\f{a}{b},$$
which is always true by (3.1).
If $a<b$, then by (3.1) we have that 
$$\f{\vartheta_1(a)}{a}\le\f{\vartheta_2(b)}{b}\,\,\Leftrightarrow\,\,\f{b}{a}\le\f{\vartheta_2(b)}{\vartheta_1(a)},$$
which implies the inequality (3.2). Hence we are done. \qed

\begin{lemma}$($\cite{K4}$)$ $(a)$ If $\ap,\bt\in\BR$, $A,B\ge 0$ and $p>1$, then we have that
\begin{equation*}|\bt-\ap|^{p-2}(\bt-\ap)(\bt B^p-\ap A^p)\ge-c_p\,\bigl(|\ap|+|\bt|\bigr)^p\,|B-A|^p,
\end{equation*} where $c_p=\f{1}{2}(p-1)^{p-1}$.

$(b)$ If $\ap,\bt\in\BR$ with $\bt\ge\ap$, $A,B\ge 0$ and $p>1$, then we have that
\begin{equation*}(\bt-\ap)^{p-1}(\bt B^p-\ap A^p)\ge\f{1}{4}(\bt-\ap)^p(A^p+B^p)-d_p\,\bigl(|\ap|+|\bt|\bigr)^p\,|B-A|^p,
\end{equation*} where $d_p=\f{1}{2}[2(p-1)]^{p-1}$.

$(c)$ If $A\ge B\ge 0$ and $p>1$, then we have that  
$$(A-B)^{p-1}\ge b_p A^{p-1}-B^{p-1},$$ where 
$b_p=\mathbbm{1}_{(1,2]}(p)+2^{-(p-1)}\mathbbm{1}_{(2,\iy)}(p)$.
\end{lemma}

\,\,\,In order to prove our results, we need two lemmas \cite{GT} which are useful in applying the De Giorgi method and obtaining nonlocal weak Harnack inequalities.

\begin{lemma} Let $\{N_k\}_{k=0}^{\iy}\subset\BR$ be a sequence of positive numbers such that
$$N_{k+1}\le d_0\,e_0^k N_k^{1+\e}$$
where $d_0,\e>0$ and $e_0>1$. If $\,N_0\le d_0^{-1/\e} e_0^{-1/\e^2}$, then we have that 
$$N_k\le e_0^{-k/\e}\,N_0$$ for any $k=0,1,\cdots$ and moreover $\lim_{k\to\iy}N_k=0$.
\end{lemma}

\begin{lemma} Let $f$ be a nonnegative bounded function defined in $[t_0,t_1]$, where $0\le t_0<t_1$. Suppose that there are nonnegative constants $c_1,c_2,\theta,$ and $\e\in(0,1)$ such that
$$f(t)\le \f{c_1}{(\tau-t)^{\theta}}+c_2+\e\,f(\tau)$$ for any $t,\tau\in[t_0,t_1]$ with $t<\tau$. Then there exists a constant $c>0$ $($depending only on $\theta$ and $\e$\,$)$ such that
$$f(\rho)\le c\,\biggl[\f{c_1}{(R-\rho)^{\theta}}+c_2\biggr]$$ for any $\rho,R\in[t_0,t_1]$ with $\rho<R$.
\end{lemma}

\begin{lemma}$($\cite{DKP1}$)$ If $p\ge 1$, $\vep\in(0,1]$ and $a,b\ge 0$, then the inequality holds
$$a^p\le b^p+c_p\vep\,b^p+(1+c_p\vep)\vep^{1-p}|a-b|^p,$$
where $c_p=(p-1)\Gm(1\vee(p-2))$ for the standard Gamma function $\Gm$.
\end{lemma}

\section{ Embedding properties of the fractional Sobolev spaces }

We furnish several embedding properties of the fractional Sobolev spaces which are very useful in oncoming estimates. The first one of them is the following version of the fractional Sobolev inequality (2.3), which can be derived by using the change of variables from the fractional Sobolev inequality (see \cite{K4}).

\begin{prop} Let $B_R$ be a ball with radius $R>0$. If $\,s\in(0,1)$ and $p\in[1,\iy)$ satisfy $sp<n$, then there is a constant $c=c(n,p,s)>0$ such that
\begin{equation*}\|f\|_{L^\tau(B_R)}\le c\,R^{-n(\f{1}{p}-\f{1}{\tau})}\|f\|_{L^p(B_R)}+c\,R^{-n(\f{1}{p}-\f{1}{\tau})+s}[f]_{W^{s,p}(B_R)}
\end{equation*} for any $\tau\in[p,p_s^*]$.
In particular, if $\tau=p_s^*$, then we have that
\begin{equation}\|f\|_{L^{p_s^*}(B_R)}\le c\,R^{-s}\|f\|_{L^p(B_R)}+c\,[f]_{W^{s,p}(B_R)}.
\end{equation}
\end{prop}

\,

We next provide the Sobolev-Poincar\'e inequality as follows (see \cite{DPV,N}).

\begin{lemma} If $u\in W^{s,p}(B)$ for a ball $B\subset\BR^n$ and $1\le\tau\le p_s^*$, then there exists a constant $c=c(n,p,s,\tau)>0$ such that
\begin{equation}\biggl(\,\intavg_B|u-u_B|^\tau\,dx\biggr)^{\f{1}{\tau}}\le c\,|B|^{\f{s}{n}}\biggl(\,\intavg_B\int_B\f{|u(x)-u(y)|^p}{|x-y|^{n+ps}}\,dx\,dy\biggr)^{\f{1}{p}}.
\end{equation}
\end{lemma}

\section{ Nonlocal Caccioppoli type inequality }
In this section, we obtain nonlocal Caccioppoli type inequalities which is one of key steps for the proof of nonlocal Harnack inequality of weak solutions to the nonlocal double phase equation (1.5). This is also a very useful tool in obtaining the local boundedness and an interior H\"older continuity of weak solutions to the nonlocal double phase equation (1.5).

\begin{lemma} Let $qt\le ps\,$ for $s,t\in(0,1)$ and $1<p\le q<\iy$, $g\in L^{q-1}_{qt}(\BR^n)$ and $B^0_{2r}\subset\Om$. If $u\in\fD_g^-(\Om)$ is a weak subsolution to the nonlocal double phase equation $(1.5)$, then for any nonnegative $\vp\in C^\iy_c(B^0_r)$ we have the estimate
\begin{equation*}\begin{split}
&\iint_{B^0_r\times B^0_r}\biggl(|w(x)\phi(x)-w(y)\phi(y)|^p K_{ps}(x,y) \\
&\quad\qquad\qquad\qquad+\fa(x,y)|w(x)\vp(x)-w(y)\vp(y)|^q K_{qt}(x,y)\biggr)\,dx\,dy \\
&\quad\ls\iint_{B^0_r\times B^0_r}\biggl((w^p(x)+w^p(y))|\phi(x)-\phi(y)|^p K_{ps}(x,y) \\
&\quad\qquad\qquad\qquad+\fa(x,y)(w^q(x)+w^q(y))|\vp(x)-\vp(y)|^q K_{qt}(x,y)\biggr)\,dx\,dy \\
&\,\,\quad+\biggl(\,\sup_{x\in\supp(\vp)}\int_{\BR^n\s B^0_r}\bigl[w^{p-1}(y) K_{ps}(x,y)+\fa(x,y)w^{q-1}(y) K_{qt}(x,y)\bigr]\,dy\biggr)\int_{B^0_r}w(x)\vp^q(x)\,dx
\end{split}\end{equation*} where $\phi=\vp^{q/p}$ and $w=(u-M)_+$ for $M\in(0,\iy)$.
\end{lemma}

\pf For simplicity, we assume that $x_0=0$. Let $w=(u-M)_+$ for $M\in[0,\iy)$ and take any nonnegative $\vp\in C^\iy_c(B_r)$. We now use 
$\zt=w\vp^q$ as a testing function in the weak formulation of the nonlocal equation.
Then we have that
\begin{equation}\begin{split}
\la u,\zt\ra_\cL&=\iint_{\BR^n\times\BR^n}\biggl(H_p(u(x)-u(y))(\zt(x)-\zt(y)) K_{ps}(x,y) \\
&\qquad\qquad+\fa(x,y)H_q(u(x)-u(y))(\zt(x)-\zt(y)) K_{qt}(x,y)\biggr)\,dx\,dy\le 0.
\end{split}\end{equation}
The above inequality can be splitted into two parts as follows;
\begin{equation}\begin{split}
\la u,\zt\ra_\cL&=\iint_{B_r\times B_r}\biggl(H_p(u(x)-u(y))(\zt(x)-\zt(y)) K_{ps}(x,y) \\
&\qquad\qquad+\fa(x,y)H_q(u(x)-u(y))(\zt(x)-\zt(y)) K_{qt}(x,y)\biggr)\,dx\,dy \\
&\quad+2\iint_{(\BR^n\s B_r)\times B_r}\biggl(H_p(u(x)-u(y))\zt(x) K_{ps}(x,y) \\
&\quad\qquad\qquad\qquad+\fa(x,y)H_q(u(x)-u(y))\zt(x) K_{qt}(x,y)\biggr)\,dx\,dy \\
&:=I_1+2\, I_2.
\end{split}\end{equation}
For the estimate of $I_1$, without loss of generality we may assume that $u(x)\ge u(y).$
Then we see that $w(x)\ge w(y)$ and 
\begin{equation}H_h(u(x)-u(y))(\zt(x)-\zt(y))\ge(w(x)-w(y))^{h-1}(\zt(x)-\zt(y))
\end{equation} for $h\in\{p,q\}$, whenever $x,y\in B_r$; indeed, it can immediately be checked by considering three possible cases (i) $u(x),u(y)>M$, (ii) $u(x)>M$, $u(y)\le M$ and (iii) $u(y)\le u(x)\le M$. 
We also note that
\begin{equation}\begin{split}
|\phi(x)w(x)-\phi(y)w(y)|^p&\le 2^{p-1}(w(x)-w(y))^p(\phi^p(x)+\phi^p(y)) \\
&\quad+2^{p-1}(w^p(x)+w^p(y))|\phi(x)-\phi(y)|^p
\end{split}\end{equation}
and 
\begin{equation}\begin{split}
|\vp(x)w(x)-\vp(y)w(y)|^q&\le 2^{q-1}(w(x)-w(y))^q(\vp^q(x)+\vp^q(y)) \\
&\quad+2^{q-1}(w^q(x)+w^q(y))|\vp(x)-\vp(y)|^q.
\end{split}\end{equation}
By Lemma 3.2, (5.3), (5.4) and (5.5), we obtain that
\begin{equation}\begin{split}
H_p(u(x)-u(y))(\zt(x)-\zt(y))&\ge\f{1}{4}(w(x)-w(y))^p(\phi^p(x)+\phi^p(y)) \\
&\quad-d_p(w^p(x)+w^p(y))|\phi(x)-\phi(y)|^p \\
&\ge\f{1}{2^{p+1}}|\phi(x)w(x)-\phi(y)w(y)|^p \\
&\quad-\bigl(\f{1}{4}+d_p\bigr)(w^p(x)+w^p(y))|\phi(x)-\phi(y)|^p
\end{split}\end{equation}
and 
\begin{equation}\begin{split}
H_q(u(x)-u(y))(\zt(x)-\zt(y))&\ge\f{1}{4}(w(x)-w(y))^q(\vp^q(x)+\vp^q(y)) \\
&\quad-d_q(w^q(x)+w^q(y))|\vp(x)-\vp(y)|^q \\
&\ge\f{1}{2^{q+1}}|\vp(x)w(x)-\vp(y)w(y)|^q \\
&\quad-\bigl(\f{1}{4}+d_q\bigr)(w^q(x)+w^q(y))|\vp(x)-\vp(y)|^q.
\end{split}\end{equation}
Thus it follows from (5.6) and (5.7) that
\begin{equation}\begin{split}
I_1&\ge\f{1}{2^{q+1}}\iint_{B_r\times B_r}\biggl(|\phi(x)w(x)-\phi(y)w(y)|^p K_{ps}(x,y)  \\
&\quad\qquad\qquad\qquad+\fa(x,y)|\vp(x)w(x)-\vp(y)w(y)|^q K_{qt}(x,y)
\biggr)\,dx\,dy \\
&\quad -2^q\bigl(\f{1}{4}+d_q\bigr)\iint_{B_r\times B_r}\biggl((w^p(x)+w^p(y))|\phi(x)-\phi(y)|^p \\
&\qquad\qquad\qquad\qquad\qquad+\fa(x,y)(w^q(x)+w^q(y))|\vp(x)-\vp(y)|^q\biggr)\,dx\,dy.
\end{split}\end{equation}
For the estimate of $I_2$, we observe that
\begin{equation*}\begin{split}
H_p(u(x)-u(y))\zt(x)&\ge-(u(y)-u(x)_+^{p-1}(u(x)-M)_+\phi^p(x) \\
&\ge-(u(y)-M)_+^{p-1}(u(x)-M)_+\phi^p(x) \\
&=-w^{p-1}(y)w(x)\phi^p(x)
\end{split}\end{equation*}
and 
\begin{equation*}\begin{split}
H_q(u(x)-u(y))\zt(x)&\ge-(u(y)-u(x)_+^{q-1}(u(x)-M)_+\vp^q(x) \\
&\ge-(u(y)-M)_+^{q-1}(u(x)-M)_+\vp^q(x) \\
&=-w^{q-1}(y)w(x)\vp^q(x).
\end{split}\end{equation*}
Thus we have that
\begin{equation}\begin{split}
I_2&\ge-\iint_{(\BR^n\s B_r)\times B_r}\biggl(w^{p-1}(y)K_{ps}(x,y)+\fa(x,y)w^{q-1}(y)K_{qt}(x,y)\biggr)\,w(x)\vp^q(x)\,dx\,dy \\
&\ge-\biggl(\sup_{x\in Q}\int_{\BR^n\s B_r}\bigl[w^{p-1}(y)K_{ps}(x,y)+\fa(x,y)w^{q-1}K_{qt}(x,y)\bigr]\,dy\biggr)\int_{B_r}w(x)\vp^q(x)\,dx
\end{split}\end{equation} where $Q=\supp(\vp)$.
Therefore the required inequality can be derived from (5.1), (5.2), (5.8) and (5.9). 
Hence we complete the proof. \qed

\begin{lemma} Let $qt\le ps\,$ for $s,t\in(0,1)$ and $1<p\le q<\iy$, $g\in L^{q-1}_{qt}(\BR^n)$ and $B^0_{2r}\subset\Om$. If $u\in\fD_g^+(\Om)$ is a weak supersolution to the nonlocal double phase equation $(1.5)$, then for any nonnegative $\vp\in C^\iy_c(B^0_r)$ we have the estimate
\begin{equation*}\begin{split}
&\iint_{B^0_r\times B^0_r}\biggl(|w(x)\phi(x)-w(y)\phi(y)|^p K_{ps}(x,y) \\
&\quad\qquad\qquad\qquad+\fa(x,y)|w(x)\vp(x)-w(y)\vp(y)|^q K_{qt}(x,y)\biggr)\,dx\,dy \\
&\quad\ls\iint_{B^0_r\times B^0_r}\biggl((w^p(x)+ w^p(y))|\phi(x)-\phi(y)|^p K_{ps}(x,y) \\
&\quad\qquad\qquad\qquad+\fa(x,y)(w^q(x)+ w^q(y))|\vp(x)-\vp(y)|^q K_{qt}(x,y)\biggr)\,dx\,dy \\
&\,\,\quad+\biggl(\,\sup_{x\in\supp(\vp)}\int_{\BR^n\s B^0_r}\bigl[w^{p-1}(y) K_{ps}(x,y)+\fa(x,y)w^{q-1}(y) K_{qt}(x,y)\bigr]\,dy\biggr)\int_{B^0_r}w(x)\vp^q(x)\,dx
\end{split}\end{equation*} where $\phi=\vp^{q/p}$ and $w=(u-M)_-=(M-u)_+$ for $M\in(0,\iy)$.
\end{lemma}

\pf It can be done in a similar way to the proof of Lemma 5.1. \qed

\,

Next, we will obtain the local boundedness of weak subsolutions to the nonlocal double phase equations and a relation between the nolocal tail terms of the positive part and the negative part of its weak subsolutions by using the De Giorgi-Nash-Moser theory.

\,

{\bf [Proof of Theorem 1.3]} (a) First, we deal with the case $qt\le ps<n$ and 
$$1<p\le q\le\f{np}{n-ps}.$$
Take any $\vp\in C^\iy_c(B^0_r)$ with $|\n\vp|\ls 1/r$ on $\BR^n$. Let $w=(u-M)_+$ for $M\in(0,\iy)$. Set $\phi=\vp^{q/p}$. Then we see that 
\begin{equation}\n\phi=\f{q}{p}\,\vp^{\f{q}{p}-1}\n\vp.
\end{equation}
Since $0<\ei\fa\le\|\fa\|_\iy$, we see that $\|\fa\|_\iy\neq 0$. Using symmetry, (5.10), Lemma 5.1 and the mean value theorem, we have that
\begin{equation}\begin{split}
&\iint_{B^0_r\times B^0_r}\biggl(|w(x)\phi(x)-w(y)\phi(y)|^p K_{ps}(x,y) \\
&\qquad+\fa(x,y)|w(x)\vp(x)-w(y)\vp(y)|^q K_{qt}(x,y)\biggr)\,dx\,dy  \\
&\ls r^{p-ps}\|\n\phi\|^p_{L^\iy(B^0_r)}\|w\|^p_{L^p(B^0_r)}+\|\fa\|_\iy\,r^{p-pt}\|\n\vp\|^q_{L^\iy(B^0_r)}\|w\|^q_{L^q(B^0_r)} \\
&\qquad+\cA(w,\vp;r,s,t)\,\|w\|_{L^1(B^0_r)}
\end{split}\end{equation} where $\cA(w,\vp;r,s,t)$ is the value given by
$$\cA(w,\vp;r,s,t)=\sup_{x\in\supp(\vp)}\int_{\BR^n\s B^0_r}\bigl[w^{p-1}(y) K_{ps}(x,y)+\|\fa\|_\iy\,w^{q-1}(y) K_{qt}(x,y)\bigr]\,dy.$$
By (1.4), we have that 
\begin{equation}\begin{split}\iint_{B^0_r\times B^0_r}&\biggl(|w(x)\phi(x)-w(y)\phi(y)|^p K_{ps}(x,y) \\
&\qquad+\fa(x,y)|w(x)\vp(x)-w(y)\vp(y)|^q K_{qt}(x,y)\biggr)\,dx\,dy  \\
&\ge\|\fa\|_* \iint_{B^0_r\times B^0_r}\biggl(|w(x)\phi(x)-w(y)\phi(y)|^p K_{ps}(x,y) \\
&\quad\qquad+\|\fa\|_\iy\,|w(x)\vp(x)-w(y)\vp(y)|^q K_{qt}(x,y)\biggr)\,dx\,dy,
\end{split}\end{equation} because $0<\|\fa\|_*:=\ei\fa/\|\fa\|_\iy\le 1$.
Thus, by (5.11) and (5.12), we have the estimate
\begin{equation}\begin{split}
\iint_{B^0_r\times B^0_r}&\biggl(|w(x)\phi(x)-w(y)\phi(y)|^p K_{ps}(x,y) \\
&\qquad+\|\fa\|_\iy\,|w(x)\vp(x)-w(y)\vp(y)|^q K_{qt}(x,y)\biggr)\,dx\,dy  \\
&\ls r^{p-ps}\|\n\phi\|^p_{L^\iy(B^0_r)}\|w\|^p_{L^p(B^0_r)}+\|\fa\|_\iy\,r^{q-qt}\|\n\vp\|^q_{L^\iy(B^0_r)}\|w\|^q_{L^q(B^0_r)} \\
&\qquad+\cA(w,\vp;r,s,t)\,\|w\|_{L^1(B^0_r)}
\end{split}\end{equation} where $\cA(w,\vp;r,s,t)$ is the value given by
$$\cA(w,\vp;r,s,t)=\sup_{x\in\supp(\vp)}\int_{\BR^n\s B^0_r}\bigl[w^{p-1}(y) K_{ps}(x,y)+\|\fa\|_\iy\,w^{q-1}(y) K_{qt}(x,y)\bigr]\,dy.$$
From (5.13), we have that
\begin{equation}\begin{split}
&\iint_{B^0_r\times B^0_r}\biggl(\,\f{|w(x)\phi(x)-w(y)\phi(y)|^p}{|x-y|^{n+ps}}+\|\fa\|_\iy\,\f{|w(x)\vp(x)-w(y)\vp(y)|^q}{|x-y|^{n+qt}}\biggr)\,dx\,dy \\
&\ls r^{p-ps}\|\n\phi\|^p_{L^\iy(B^0_r)}\|w\|^p_{L^p(B^0_r)}+\|\fa\|_\iy\,r^{q-qt}\|\n\vp\|^q_{L^\iy(B^0_r)}\|w\|^q_{L^q(B^0_r)} \\
&\qquad+\cA(w,\vp;r,s,t)\,\|w\|_{L^1(B^0_r)}.
\end{split}\end{equation}
Applying Proposition 4.1 to (5.14), we obtain that
\begin{equation}\begin{split}
\biggl(\,\intavg_{B^0_r}|w\phi|^{p\gm}\,dx\biggr)^{\f{1}{\gm}}&\ls \bigl(r^{p-ps}\|\n\phi\|^p_{L^\iy(B^0_r)}+r^{-ps}\bigr)r^{ps}\intavg_{B^0_r}|w|^p\,dx \\
&\quad+r^{q-qt}\|\n\vp\|^q_{L^\iy(B^0_r)}\,r^{ps}\,\|\fa\|_\iy\intavg_{B^0_r}|w|^q\,dx \\
&\quad+\cA(w,\vp;r,s,t)r^{ps}\intavg_{B^0_r}w\,dx
\end{split}\end{equation}
where the constant $\gm$ is given by
\begin{equation*}\gm:=\f{n}{n-ps}>1.\end{equation*}
For $k=0,1,2,\cdots$, we set
\begin{equation*}\begin{split}r_k&=(1+2^{-k})r, \,\,r_k^*=\f{r_k+r_{k+1}}{2}, \\
M_k&=M+(1-2^{-k})M_*,\,\,M_k^*=\f{M_k+M_{k+1}}{2},
\end{split}\end{equation*}
$w_k=(u-M_k)_+$ and $w_k^*=(u-M_k^*)_+$ for a constant $M_*>0$ to be selected later, and we set
$$N_k=\intavg_{B^0_{r_k}}\Theta(w_k)\,dx.$$
In (5.15), for each $k=0,1,2,\cdots$, take a function $\vp_k\in C_c^\iy(B^0_{r^*_k})$ with $\vp_k|_{B^0_{r_{k+1}}}\equiv 1$ such that $0\le\vp_k\le 1$ and 
\begin{equation}|\n\vp_k|\ls\f{2^{k+2}}{r}\,\text{ in $\BR^n$ }\end{equation}
and we set $\phi_k=\vp^{q/p}_k$. Then by (5.10) we see that
\begin{equation}|\n\phi_k|\ls\f{2^{k+2}}{r}\,\text{ in $\BR^n$. }\end{equation}
Since $w_k\ge w_k^*\ge w_{k+1}$ and $w_k^*\ge M_{k+1}-M_k^*=2^{-k-2}M_*$ whenever $u(x)\ge M_{k+1}$, we have that
\begin{equation}\begin{split}
N_{k+1}&\ls\f{1}{|B^0_{r_k}|}\int_{B^0_{r_{k+1}}}\f{w^p_{k+1}(w^*_k)^{p(\gm-1)}}{(M_{k+1}-M^*_k)^{p(\gm-1)}}\,dx
+\f{\|\fa\|_\iy}{|B^0_{r_k}|}\int_{B^0_{r_{k+1}}}\f{w^q_{k+1}(w^*_k)^{p(\gm-q/p)}}{(M_{k+1}-M^*_k)^{p(\gm-q/p)}}\,dx \\
&\ls\biggl(\f{2^k}{M_*}\biggr)^{p(\gm-1)}\intavg_{B^0_{r_k}}|w^*_k\phi_k|^{p\gm}\,dx
+\biggl(\f{2^k}{M_*}\biggr)^{p(\gm-q/p)}\intavg_{B^0_{r_k}}|w^*_k\phi_k|^{p\gm}\,dx.
\end{split}\end{equation}
In (5.18), we used the fact that $\gm-q/p\ge 0$, i.e. $q\le\f{np}{n-ps}$.
Also we note that 
\begin{equation}\cA(w^*_k,\vp_k;r_k,s,t)\le c\, 2^{k(n+ps)}r^{-ps}\,\cT_r(w_0;x_0),
\end{equation}
because $\vp_k\in C_c^\iy(B_{r_k})$, $w^*_k\le w_0$ for all $k$ and 
$$|y-x|\ge|y-x_0|-|x-x_0|\ge\bigl(1-\f{r^*_k}{r_k}\bigr)\,|y-x_0|\ge 2^{-k-2}|y-x_0|$$
for any $x\in B_{r^*_k}$ and $y\in\BR^n\s B_{r_k}$.
Let $k_0\in\BN$ be the least natural number so that $2^k\ge M_*$ for any $k\ge k_0$. If $k\ge k_0$, then it follows from (5.15), (5.16), (5.17), (5.18) and  (5.19) that
\begin{equation}\begin{split}
\biggl(\f{2^k}{M_*}\biggr)^{-\f{p(\gm-1)}{\gm}}N_{k+1}^{\f{1}{\gm}}&\ls\biggl(\,\intavg_{B^0_{r_k}}|w^*_k\phi_k|^{p\gm}\,dx\biggr)^{\f{1}{\gm}} \\
&\le c\,2^{pk}\intavg_{B^0_{r_k}}|w^*_k|^p\,dx+c\,2^{qk}\|\fa\|_\iy\intavg_{B^0_{r_k}}|w^*_k|^q\,dx \\
&\quad+2^{k(n+ps)}\biggl(\f{r_k}{r}\biggr)^{ps}\cT_r(w_0;x_0)\intavg_{B^0_{r_k}}w^*_k\,dx \\
&\le c\,2^{qk}N_k+c\,\cT_r(w_0;x_0)\,2^{k(n+ps)}\intavg_{B^0_{r_k}}\f{w^*_k w_k^{p-1}}{(M^*_k-M_k)^{p-1}}\,dx \\
&\le c\,\biggl(\,2^{qk}+2^{k(n+ps)}\biggl(\f{2^k}{M_*}\biggr)^{p-1}\cT_r(w_0;x_0)\biggr)\,N_k,
\end{split}\end{equation}
because we see that $0\le w^*_k\le w_k$ and $w_k(x)\ge M^*_k-M_k=2^{-k-2}M_*$ if $u(x)\ge M^*_k$. Taking $M^*$ in the above so that
\begin{equation*}M_*\ge\dt\,[\cT_r(w_0;x_0)]^{\f{1}{p-1}}\,\,\text{ for $\dt\in(0,1]\,\,$ in (5.20), }
\end{equation*}
we have that
$$\f{N_{k+1}}{M_*^p}\le d_0 \,a^k\,\biggl(\f{N_k}{M_*^p}\biggr)^{1+\e}$$ for every $k\ge k_0$, 
where $d_0=c^\gm\dt^{-(p-1)\gm}>0$, $a=2^{p(\gm-1)+n+ps+p-1+q}>1$ and $\e=\gm-1>0$. 
If $N_{k_0}\le d_0^{-\f{1}{\e}}a^{-\f{1}{\e^2}}M_*^p$, then we set
\begin{equation*}M_*=\dt\,[\cT_r(w_0;x_0)]^{\f{1}{p-1}}+c_0\,\dt^{-\f{(p-1)n}{sp^2}}a^{\f{(n-ps)^2}{p^3 s^2}}N^{\f{1}{p}}_{k_0}
\end{equation*} 
where $c_0=c^{\f{n}{sp^2}}$.
By applying Lemma 3.3, we conclude that
\begin{equation*}\begin{split}
&\sup_{B^0_r}u\le M+M_* \\
&\le M+\dt\,[\cT_r(w_0;x_0)]^{\f{1}{p-1}}+c_0\,\dt^{-\f{(p-1)n}{sp^2}}a^{\f{(n-ps)^2}{p^3 s^2}}
\biggl(\,\intavg_{B^0_{r_{k_0}}}\Theta((u-M-(1-2^{-k_0})M_*)_+)\,dx\biggr)^{\f{1}{p}} \\
&\le M+\dt\,[\cT_r(w_0;x_0)]^{\f{1}{p-1}}+c_0\,\dt^{-\f{(p-1)n}{sp^2}}a^{\f{(n-ps)^2}{p^3 s^2}}
\biggl(\,\intavg_{B^0_{2r}}\Theta((u-M)_+)\,dx\biggr)^{\f{1}{p}}.
\end{split}\end{equation*}
Letting $M\downarrow 0$ in the above, we obtain that
\begin{equation*}\sup_{B^0_r}u\le\dt\,[\cT_r(u^+;x_0)]^{\f{1}{p-1}}+C_0\,\dt^{-\f{(p-1)n}{sp^2}}
\biggl(\,\intavg_{B^0_{2r}}\Theta(u^+)\,dx\biggr)^{\f{1}{p}}.\end{equation*}

(b) Next we consider the case $qt\le ps=n$ and $1<p\le q<\iy$. If $p=q$, then by (5.12) it is exactly the same as that of the fractional $p$-Laplacian case. So, without loss of generality, we assume that  $qt\le ps=n$ and $1<p<q<\iy$. 
This proof will be proceeded in a similar way with the first case (a). Let us take any $s_1\in(0,s)$. Then we see that $ps_1<n$ and 
$$w\phi-(w\phi)_{B_r^0}\in W^{s,p}(B^0_r)\subseteq W^{s_1,p}(B^0_r),$$ where $$(w\phi)_{B^0_r}=\intavg_{B^0_r}w(x)\phi(x)\,dx.$$
For a moment, we consider
$$1<p<q<\f{np}{n-ps_1}:=p_1^*.$$
Since the function 
$$B_f(p):=\biggl(\,\intavg_E |f(x)|^p\,dx\biggr)^{1/p}$$ is increasing on $(0,\iy)$ for $f\in \cF(\BR^n)$ and a measurable set $E\subset\BR^n$ with $0<|E|<\iy$, we have that
\begin{equation}\biggl(\,\int_{B^0_r}|w\phi-(w\phi)_{B^0_r}|^q\,dx\biggr)^{\f{1}{q}}\le|B^0_r|^{\f{p_1^*-q}{p_1^* q}}
\biggl(\,\int_{B^0_r}|w\phi-(w\phi)_{B^0_r}|^{p_1^*}\,dx\biggr)^{\f{1}{p_1^*}}.
\end{equation}
Thus it follows from (5.21) that
\begin{equation}\begin{split}
&\biggl|\,\biggl(\,\intavg_{B^0_r}|w\phi-(w\phi)_{B^0_r}|^q\,dx\biggr)^{\f{1}{q}}-|(w\phi)_{B^0_r}|\,\biggr|^p
\le c\,\biggl(\,\intavg_{B^0_r}|w\phi-(w\phi)_{B^0_r}|^q\,dx\biggr)^{\f{p}{q}} \\
&\qquad\qquad\qquad\qquad\qquad\qquad\qquad\le c\,|B^0_r|^{\f{(p_1^*-q)p}{p_1^* q}-\f{p}{q}}\,
\biggl(\,\int_{B^0_r}|w\phi-(w\phi)_{B^0_r}|^{p_1^*}\,dx\biggr)^{\f{p}{p_1^*}}.
\end{split}\end{equation}
Also, by Jensen's inequality and the fact that $|a+b|^p\le 2^p(|a|^p+|b|^p)$ for $p>1$ and $a,b\in\BR$, we have that
\begin{equation}2^p\biggl|\,\biggl(\,\intavg_{B^0_r}|w\phi-(w\phi)_{B^0_r}|^q\,dx\biggr)^{\f{1}{q}}-|(w\phi)_{B^0_r}|\,\biggr|^p\ge\biggl(\,\intavg_{B^0_r}|w\phi|^q\,dx\biggr)^{\f{p}{q}}-2^p\intavg_{B^0_r}|w\phi|^p\,dx.
\end{equation}
Applying Proposition 2.2, by (5.22) and (5.23) we have that
\begin{equation*}\begin{split}
\biggl(\,\intavg_{B^0_{r_k}}|w_k^*\phi_k|^q\,dx\biggr)^{\f{p}{q}}
&\ls|B^0_{r_k}|^{\f{(p_1^*-q)p}{p_1^* q}-\f{p}{q}}\,
\biggl(\,\int_{B^0_{r_k}}|w_k^*\phi_k-(w_k^*\phi_k)_{B^0_{r_k}}|^{p_1^*}\,dx\biggr)^{\f{p}{p_1^*}} \\
&\qquad+\intavg_{B^0_{r_k}}|w_k^*\phi_k|^p\,dx \\
&\ls r_k^{ps_1-n}\iint_{B^0_{r_k}\times B^0_{r_k}}\f{|w_k^*(x)\phi_k(x)-w_k^*(y)\phi_k(y)|^p}{|x-y|^{n+ps_1}}\,dx\,dy \\
&\qquad+\intavg_{B^0_{r_k}}|w_k^*\phi_k|^p\,dx \\ 
\end{split}\end{equation*}
for any $s_1\in(0,s)$. Since $|x-y|^{n+ps}\ls |x-y|^{n+ps_1}$ for any $x,y\in B^0_r$ and any $s_1\in(0,s)$, taking the limit $s_1\nearrow s$ just in the above it follows from Lebesgue dominated convergence theorem that
\begin{equation*}\begin{split}
\biggl(\,\intavg_{B^0_{r_k}}|w_k^*\phi_k|^q\,dx\biggr)^{\f{p}{q}} &\ls r_k^{ps-n}\iint_{B^0_{r_k}\times B^0_{r_k}}\f{|w_k^*(x)\phi_k(x)-w_k^*(y)\phi_k(y)|^p}{|x-y|^{n+ps}}\,dx\,dy \\
&\qquad+\intavg_{B^0_{r_k}}|w_k^*\phi_k|^p\,dx.
\end{split}\end{equation*}
For $k=0,1,2,\cdots$, we set $N_k=\ds\intavg_{B^0_{r_k}}\Theta(w_k)\,dx$. Let $k_0\in\BN$ be the least natural number so that $2^k\ge M_*$ for any $k\ge k_0$.
Since $ps\ge qt$, $w_k\ge w_k^*\ge w_{k+1}$ and $w_k^*\ge M_{k+1}-M_k^*=2^{-k-2}M_*$, 
by (5.14) and (5.19) we have that
\begin{equation}\begin{split}
&N_{k+1}^{\f{p}{q}}
\ls\biggl[\biggl(\f{2^k}{M_*}\biggr)^{q-p}+1\biggr]^{\f{p}{q}}\biggl(\,\intavg_{B^0_{r_k}}|w_k^*\phi_k|^q\,dx\biggr)^{\f{p}{q}} \\
&\ls\biggl[\biggl(\f{2^k}{M_*}\biggr)^{q-p}+1\biggr]^{\f{p}{q}} \biggl[(r_k^{p-ps}\|\phi_k\|^p_{L^\iy(B^0_{r_k})}\,r_k^{ps}+1)\intavg_{B^0_{r_k}}|w_k|^p\,dx \\
&\qquad+ r_k^{q-qt}\|\vp_k\|^q_{L^\iy(B^0_{r_k})}\,r_k^{ps}\intavg_{B^0_{r_k}}|w_k|^q\,dx
+r_k^{ps}\cA(w_k^*,\phi_k;r_k,s,t)\intavg_{B^0_{r_k}}w_k\,dx\biggr] \\
&\ls\biggl[\biggl(\f{2^k}{M_*}\biggr)^{q-p}+1\biggr]^{\f{p}{q}}\,\biggl[N_k
+2^{k(n+ps)}\biggl(\f{r_k}{r}\biggr)^{ps}
\cT_r(w_0;x_0)\intavg_{B^0_{r_k}}\f{w_k(w_k^*)^{p-1}}{(M_k^*-M_k)^{p-1}}\,dx\biggr] \\
&\ls\biggl[\biggl(\f{2^k}{M_*}\biggr)^{q-p}+1\biggr]^{\f{p}{q}}\,\biggl[1+2^{k(n+ps)}\biggl(\f{2^k}{M_*}\biggr)^{p-1}\cT_r(w_0;x_0)\biggr]\,N_k.
\end{split}\end{equation}
Taking $M_*$ in (5.24) so that
\begin{equation*}M_*\ge\dt\,[\cT_r(w_0;x_0)]^{\f{1}{p-1}}\,\,\text{ for $\dt\in(0,1]$, }
\end{equation*}
then we have that
\begin{equation*}N_{k+1}^{\f{p}{q}}\le c\,\biggl(\f{2^k}{M_*}\biggr)^{(q-p)\f{p}{q}}\dt^{-(p-1)}2^{k(n+ps+p-1)} N_k;
\end{equation*}
that is to say, we obtain that
\begin{equation*}\f{N_{k+1}}{M_*^p}\le b_0\,a^k\,\biggl(\f{N_k}{M_*^p}\biggr)^{1+\e}
\end{equation*} for every $k\ge k_0$,
where $b_0=c^{\f{q}{p}}\dt^{-\f{q(p-1)}{p}}>0$, $a=2^{(q-p)+\f{q(n+ps+p-1)}{p}}>1$ and $\e=q/p-1>0$.
If $N_{k_0}\le b_0^{-\f{1}{\e}}a^{-\f{1}{\e^2}}M_*^p$, then  we set
$$M_*=\dt\,[\cT_r(w_0;x_0)]^{\f{1}{p-1}}+c_0\,\dt^{-\f{q(p-1)}{p(q-p)}}\,a^{\f{p}{(q-p)^2}}\,N_{k_0}^{\f{1}{p}}$$
where $c_0=c^{\f{q}{q-p}}$.
By applying Lemma 3.3, we conclude that
\begin{equation*}\begin{split}
\sup_{B^0_r}u&\le M+\dt\,[\cT_r(w_0;x_0)]^{\f{1}{p-1}} \\
&\qquad+c_0\,\dt^{-\f{q(p-1)}{p(q-p)}}\,a^{\f{p}{(q-p)^2}}
\biggl(\,\intavg_{B^0_{r_{k_0}}}\Theta((u-M-(1-2^{-k_0})M_*)_+)\,dx\biggr)^{\f{1}{p}} \\
&\le M+\dt\,[\cT_r(w_0;x_0)]^{\f{1}{p-1}}+c_0\,\dt^{-\f{q(p-1)}{p(q-p)}}\,a^{\f{p}{(q-p)^2}}
\biggl(\,\intavg_{B^0_{2r}}\Theta((u-M)_+)\,dx\biggr)^{\f{1}{p}}.
\end{split}\end{equation*}
Letting $M\downarrow 0$ in the above, we obtain that
\begin{equation*}\sup_{B^0_r}u\le\dt\,[\cT_r(u^+;x_0)]^{\f{1}{p-1}}+C_0\,\dt^{-\f{q(p-1)}{p(q-p)}}
\biggl(\,\intavg_{B^0_{2r}}\Theta(u^+)\,dx\biggr)^{\f{1}{p}}.\end{equation*}
Hence we obtain the required inequality. \qed

\begin{lemma} Suppose that $qt\le ps\,$ for $\,s,t\in(0,1)\,$, $1<p\le q<\iy\,$ and $\,g\in L^{q-1}_{qt}(\BR^n)$. If $u\in\fD^-_g(\Om)$ is a weak subsolution of the nonlocal double phase equation $(1.5)$ such that $u\ge 0$ in $B_R^0\subset\Om$, then we have the estimate
\begin{equation*}\cT_r(u^+;x_0)\ls\biggl[\sup_{B^0_r}u\biggr]^{p-1}+\biggl[\sup_{B^0_r}u\biggr]^{q-1} 
+\biggl(\f{r}{R}\biggr)^{ps}\cT_R(u^-;x_0)
\end{equation*}
for any $r\in(0,R)$.
\end{lemma}

\pf For simplicity, we set $x_0=0$, i.e. $B_r=B^0_r$ for $r\in(0,R)$.
Let $M=\sup_{B_r}u$ and 
$$\vp(x)=w(x)\zt^q(x)$$ for $w(x)=u(x)-2M$ and a function $\zt\in C^\iy_c(B_{3r/4})$ satisfying that 
$$\zt|_{B_{r/2}}\equiv 1, \,\,0\le\zt\le 1 \,\,\text{ and $\,\,|\n\zt|\ls\f{1}{r}\,$ in $\BR^n$.}$$ Then it is easy to check that $\vp\in\fD_0(\Om)$. Thus we are going to use $\vp$ as a testing function in the weak formulation of the nonlocal equation (1.5). Let $\xi=\zt^{q/p}$.
Then we have that
\begin{equation}\begin{split}
0&\ge\iint_{B_r\times B_r}\biggl(H_p(u(x)-u(y))(\vp(x)-\vp(y)) K_{ps}(x,y) \\
&\qquad\qquad+\fa(x,y)H_q(u(x)-u(y))(\vp(x)-\vp(y)) K_{qt}(x,y)\biggr)\,dx\,dy \\
&\quad+2\iint_{(\BR^n\s B_r)\times B_r}\biggl(H_p(u(x)-u(y))(u(x)-2M)\zt^q(x) K_{ps}(x,y) \\
&\quad\qquad\qquad\qquad+\fa(x,y)H_q(u(x)-u(y))(u(x)-2M)\zt^q(x) K_{qt}(x,y)\biggr)\,dx\,dy \\
&:=\cI_1+2\, \cI_2.
\end{split}\end{equation}
Since $-2M\le w(x)=u(x)-2M\le-M$ for any $x\in B_r$, it follows from (a) of Lemma 3.2 that
\begin{equation*}
H_p(w(x)-w(y))(w(x)\xi^p(x)-w(y)\xi^p(y))\ge-c_p 4^p(\xi(x)-\xi(y))^p
\end{equation*} and
\begin{equation*}
H_q(w(x)-w(y))(w(x)\zt^q(x)-w(y)\zt^q(y))\ge-c_q 4^q(\zt(x)-\zt(y))^q
\end{equation*}
for any $x,y\in B_r$. Thus we have that
\begin{equation}\begin{split}
\cI_1&\ge-c_p 4^p M^p\iint_{B_r\times B_r}(\xi(x)-\xi(y))^p K_{ps}(x,y)\,dx\,dy \\
&\qquad-c_q 4^q M^q\|\fa\|_{L^\iy}\iint_{B_r\times B_r}(\zt(x)-\zt(y))^q K_{qt}(x,y)\,dx\,dy \\
&\gs-M^p r^{-ps}|B_r|-M^q r^{-qt}|B_r|\gs-(M^p+M^q)r^{-ps}|B_r|.
\end{split}\end{equation}
Also the lower estimate of $\cI_2$ can be decomposed into two parts as follows;
\begin{equation}\begin{split}
\cI_2&\ge 4M\iint_{(\BR^n\s B_r)\times B_r}\bigl[(u(y)-M)_+^{p-1}\zt^q(x)K_{ps}(x,y)  \\
&\qquad\qquad\qquad\qquad+\fa(x,y)(u(y)-M)_+^{q-1}\zt^q(x)K_{qt}(x,y)\bigr]\,dx\,dy \\
&\quad-4M\iint_{E_M\times B_r}\bigl[(u(x)-u(y))_+^{p-1}\zt^q(x)K_{ps}(x,y) \\
&\qquad\qquad\qquad\qquad+\fa(x,y)(u(x)-u(y))_+^{q-1}\zt^q(x)K_{qt}(x,y)\bigr]\,dx\,dy \\
&:=\cI_{2,1}-\cI_{2,2}
\end{split}\end{equation}
where $E_M=\{y\in\BR^n\s B_r:u(y)<M\}$. Since we see that $$(u(y)-M)_+\ge u^+(y)-M$$ for any $y\in\BR^n\s B_r,\,$ by (c) of Lemma 3.2 we have that
\begin{equation}(u(y)-M)^{h-1}_+\ge b_h[u^+(y)]^{h-1}-M^{h-1}\end{equation}
where the constant $b_h$ is given by
$$b_h=\mathbbm{1}_{(1,2]}(h)+2^{-(h-1)}\mathbbm{1}_{(2,\iy)}(h)\,\,\text{ for $h\in\{p,q\}$.}$$ 
So by (5.28) the lower estimate of $\cI_{2,1}$ can be done as
\begin{equation}\begin{split}
\cI_{2,1}&\ge d_2\, M\,|B_r|\,I(r)-d_3 (M^p r^{-ps}+M^q r^{-qt})|B_r|
\end{split}\end{equation}
with universal constants $d_2>0$, where $I(r)$ is given by
$$I(r):=\int_{\BR^n\s B_r}\bigl(\f{[u^+(y)]^{p-1}}{|y-x_0|^{n+ps}}+\intavg_{B_{r/2}}\fa(x,y)\,dx\,\f{[u^+(y)]^{q-1}}{|y-x_0|^{n+qt}}\bigr)\,dy.$$
We note that, for any $r\in(0,R)$,
\begin{equation}\begin{split}
I(r)&\le\int_{\BR^n\s B_r}\bigl(\f{[u^+(y)]^{p-1}}{|y-x_0|^{n+ps}}+\|\fa\|_\iy\,\f{[u^+(y)]^{q-1}}{|y-x_0|^{n+qt}}\bigr)\,dy \\
&\ls\int_{\BR^n\s B_r}\f{[u^+(y)]^{q-1}+1}{|y-x_0|^{n+qt}}\,dy<\iy.
\end{split}\end{equation}
Then, by (1.4), we have that 
\begin{equation}I(r)\ge\|\fa\|_*\int_{\BR^n\s B_r}\bigl(\f{[u^+(y)]^{p-1}}{|y-x_0|^{n+ps}}+\|\fa\|_\iy\,\f{[u^+(y)]^{q-1}}{|y-x_0|^{n+qt}}\bigr)\,dy.
\end{equation}
Then it follows from (5.29) and (5.31) that
\begin{equation}\begin{split}
\cI_{2,1}&\gs d_2 M\,|B_r|\int_{\BR^n\s B_r}\bigl(\f{[u^+(y)]^{p-1}}{|y-x_0|^{n+ps}}+\|\fa\|_\iy\f{[u^+(y)]^{q-1}}{|y-x_0|^{n+qt}}\bigr)\,dy  \\
&\qquad-d_3 (M^p r^{-ps}+M^q r^{-qt})|B_r| \\
&\gs M r^{-ps}|B_r|\,\cT_r(u^+;0)-(M^p+M^q)r^{-ps} |B_r|.
\end{split}\end{equation}
If $x\in B_r$ and $y\in E_M$, then we note that, for $h\in\{p,q\}$,
\begin{equation}\begin{split}
(u(x)-u(y))_+^{h-1}&\le a_h\bigl(|u(x)-M)^{h-1}+|M-u(y)|^{h-1}\bigr)   \\
&\le a_h M^{h-1}+a_h(M+u^-(y)-u^+(y))^{h-1} \\
&\le a_h M^{h-1}+a_h(M+u^-(y))^{h-1} \\
&\le a_h(1+a_h) M^{h-1}+a_h^2[u^-(y)]^{h-1}
\end{split}\end{equation}
where $a_h=\mathbbm{1}_{(1,2]}(h)+2^{h-1}\mathbbm{1}_{(2,\iy)}(h)$, because $u^+(y)<M+u^-(y)$ for any $y\in E_M$.
Since we see that $u^-(y)=0$ for any $y\in B_R$, the upper estimate of $\cI_{2,2}$ can be obtained by
\begin{equation}\begin{split}
\cI_{2,2}&\le 4 a_p(1+a_p)M^p\iint_{(\BR^n\s B_r)\times B_r}\zt^q(x)K_{ps}(x,y)\,dx\,dy \\
&\quad+4 a_q(1+a_q)M^q\,\|\fa\|_{\iy}\iint_{(\BR^n\s B_r)\times B_r}\zt^q(x)K_{qt}(x,y)\,dx\,dy \\
&\quad+4(a_p^2\vee a_q^2)M\iint_{(\BR^n\s B_R)\times B_r}\bigl([u^-(y)]^{p-1}\zt^q(x)K_{ps}(x,y) \\
&\qquad\qquad\qquad\qquad\qquad+\fa(x,y)[u^-(y)]^{q-1}\zt^q(x)K_{qt}(x,y) \bigr)\,dx\,dy \\
&\le d_3 (M^p r^{-ps}+\|\fa\|_\iy\,M^q r^{-qt})|B_r| \\
&\quad\qquad\qquad+d_4 M|B_r|\int_{\BR^n\s B_R}\bigl(\f{[u^-(y)]^{p-1}}{|y-x_0|^{n+ps}}+\|\fa\|_\iy\f{[u^-(y)]^{q-1}}{|y-x_0|^{n+qt}}\bigr)\,dy  \\
&\ls (M^p+M^q) r^{-ps}|B_r|+M R^{-ps}|B_r|\,\cT_R(u^-;x_0)\\
\end{split}\end{equation}
Therefore the required estimate can be achieved by (5.26), (5.27), (5.29), (5.32) and (5.34). \qed

\begin{thm} Let $qt\le ps\le n\,$ for $s,t\in(0,1)$ and $1<p\le q\le p_s^*$, $g\in L^{q-1}_{qt}(\BR^n)$.  If $u\in\fD_g^-(\Om)$ is a globally nonnegative weak subsolution to the nonlocal double phase equation $(1.5)$, then there is a constant $C_0>0$ depending only on $n,s,p,q,\ld,\Ld,\Om$ and $\fa$ such that
$$\sup_{B^0_r}u\le C\,\biggl(\,\intavg_{B^0_{2r}}\Theta(u)\,dx\biggr)^{\f{1}{p}}$$ 
for any $r>0$ with $B^0_{2r}\subset\Om$,
where the function $\Theta$ is defined by
\begin{equation*}\Theta(z):=z^p+\|a\|_{\iy}\,z^q,\,\,z\ge 0.
\end{equation*}
\end{thm} 

\pf Choose some small $\dt\in(0,1]$ so that $1-\dt\, ((2 d_0)^{\f{1}{p-1}}a_p\vee (2 d_0)^{\f{1}{q-1}}a_q)>0$ and take any $r>0$ with $B^0_{2r}\subset\Om$,
where $d_0>0$ is the universal constant in Theorem 5.3 and
$$a_p=\mathbbm{1}_{[2,\iy)}(p)+2^{\f{1}{p-1}}\mathbbm{1}_{(1,2)}(p).$$
We consider two possible cases, i.e. $\sup_{B^0_r}u\le 1$ or $\,\sup_{B^0_r}u>1$.
If $\sup_{B^0_r}u\le 1$, then it follows from Theorem 1.3 and Theorem 5.3 that
\begin{equation*}\begin{split}
\sup_{B^0_r}u&\le\dt\,(2 d_0)^{\f{1}{p-1}}\,a_p\biggl[\,\sup_{B^0_r}u+2^{-\f{ps}{p-1}}\cT_{2r}(u^-;x_0)\biggr]
+C_0\,\dt^{-\f{(p-1)n}{sp^2}}\biggl(\,\intavg_{B^0_{2r}}\Theta(u)\,dx\biggr)^{\f{1}{p}}.
\end{split}\end{equation*}
Since $\cT_{2r}(u^-;x_0)=0$, the required inequality can easily derived by taking
$$C=\f{C_0\,\dt^{-\f{(p-1)n}{sp^2}}}{1-\dt\, ((2 d_0)^{\f{1}{p-1}}a_p\vee (2 d_0)^{\f{1}{q-1}}a_q)}.$$ 
If $\sup_{B^0_r}u> 1$, then it follows from Theorem 1.3 and Theorem 5.3 that
\begin{equation*}\begin{split}
\sup_{B^0_r}u&\le\dt\,(2 d_0)^{\f{1}{q-1}}\,a_q\biggl[\,\sup_{B^0_r}u+2^{-\f{ps}{q-1}}\cT_{2r}(u^-;x_0)\biggr] +C_0\,\dt^{-\f{(p-1)n}{sp^2}}\biggl(\,\intavg_{B^0_{2r}}\Theta(u)\,dx\biggr)^{\f{1}{p}}.
\end{split}\end{equation*}
Since $\cT_{2r}(u^-;x_0)=0$, the required inequality can easily derived by taking
$$C=\f{C_0\,\dt^{-\f{(p-1)n}{sp^2}}}{1-\dt\, ((2 d_0)^{\f{1}{p-1}}a_p\vee (2 d_0)^{\f{1}{q-1}}a_q)}.$$ 
Hence we complete the proof. \qed

\section{ The logarithm of a weak solution is a locally BMO function }

In this section, we shall prove that the logarithm of a weak supersolution $u$ to the nonlocal double phase equation (1.5) with $u\ge 0$ in $B^0_R\subset\Om$ is a function with locally bounded mean oscillation. To realize it, the following tool which is called the {\it fractional Poincar\'e inequality} is very useful.

Let $n\ge 1$, $p>1$, $s\in(0,1)$ and $ps<n$. For a ball $B\subset\BR^n$, let $u_B$ denote the average of $u\in W^{s,p}(B)$ over $B$, i.e.
$$u_B=\intavg_B u(y)\,dy.$$
Then it is well-knowm in \cite{BBM,MS} that 
\begin{equation}\|u-u_B\|_{L^p(B)}\le\f{c_{n,p}(1-s)|B|^{\f{s}{n}}}{(n-ps)^{(p-1)/p}}\,[u]_{W^{s,p}(B)}
\end{equation}
where $c_{n,p}>0$ is a universal constant depending only on $n$ and $p$. This inequality is very useful in obtaining the logarithmic estimate of its weak supersolutions. 

For $p>0$, we now consider a kind of local BMO spaces on $B^0_R\subset\BR^n$, whose norm $\|\cdot\|_{\bmo^p(B^0_R)}$ is defined by
$$\|u\|_{\bmo^p(B^0_R)}=\sup_{r\in(0,R/2)}\biggl(\intavg_{B^0_r}|u(y)-u_{B^0_r}|^p\,dy\biggr)^{\f{1}{p}}$$
and we write the space 
$$\bmo^p(B^0_R)=\{u\in L^1_{\loc}(\BR^n):\|u\|_{\bmo^p(B^0_R)}<\iy\}.$$ 
When $p=1$, we write $\bmo^p(B^0_R)=\bmo(B^0_R)$. Then it is easy to check that
$$\|\,|u|\,\|_{\bmo(B^0_R)}\le\|u\|_{\bmo(B^0_R)}$$
because we see that 
$$||u|-|u|_B|\le|u-u_B|.$$ Also we note that
\begin{equation}\|u\pm v\|_{\bmo(B^0_R)}\le\|u\|_{\bmo(B^0_R)}+\|v\|_{\bmo(B^0_R)}.
\end{equation}
From the fact that
$$a\wedge b=\f{a+b-|a-b|}{2}\,\,\text{ and }\,\,a\vee b=\f{a+b+|a-b|}{2}$$ for any $a,b>0$, we can easily derive that
\begin{equation}\begin{split}
\|u\wedge v\|_{\bmo(B^0_R)}&\le\|u\|_{\bmo(B^0_R)}+\|v\|_{\bmo(B^0_R)}, \\
\|u\vee v\|_{\bmo(B^0_R)}&\le\|u\|_{\bmo(B^0_R)}+\|v\|_{\bmo(B^0_R)}.
\end{split}\end{equation}
In addition, we can derive the following {\it John-Nirenberg inequality} (as in \cite{Gr}) by applying the Calder\'on-Zygmund decomposition technique in standard harmonic analysis as follows; there exists some constants $d_1,d_2>0$ depending only on the dimension $n$ such that
$$\bigl|\{x\in B^0_r:|f(x)-f_{B^0_r}|>\ld\}\bigr|\le d_1\,e^{-(d_2/\|f\|_{\bmo(B^0_R)})\ld}|B^0_r|
$$
for any $f\in\bmo(B^0_r)$, every $r>0$ with $B^0_{2r}\subset B^0_R\subset\Om$ and every $\ld>0$.
By standard analysis, this inequality makes it possible to obtain the following fact;
\begin{equation}\begin{split}
&\text{ If $f\in\bmo(B^0_R)$ for $B^0_R\subset\Om$ and $1<p<\iy$, then $\|\cdot\|_{\bmo(B^0_R)}$ is }  \\
&\,\,\,\text{ norm-equivalent to $\|\cdot\|_{\bmo^p(B^0_R)}$. }
\end{split}\end{equation}

For the next lemma, we consider various functions as follows; for $x,y\in\BR^n$, $r>0$ and $z\ge 0$, we define
\begin{equation}\begin{split}
&\Theta_*(x,y;z)=\f{z^p}{|x-y|^{ps}}+a(x,y)\f{z^q}{|x-y|^{qt}}, \,\Theta(x,y;z)=\f{z^p}{|x-y|^{ps}}+\|a\|_\iy\f{z^q}{|x-y|^{qt}}, \\
&\theta_*(x,y;z)=\f{z^{p-1}}{|x-y|^{ps}}+a(x,y)\f{z^{q-1}}{|x-y|^{qt}}, \,\theta(x,y;z)=\f{z^{p-1}}{|x-y|^{ps}}+\|a\|_\iy\f{z^{q-1}}{|x-y|^{qt}}, \\
&\Theta^*_r(z)=\f{z^p}{r^{ps}}+a(x,y) \,\f{z^q}{r^{qt}},\,\,\,\,\theta^*_r(z)=\f{z^{p-1}}{r^{ps}}+a(x,y) \,\f{z^{q-1}}{r^{qt}} ,\\
&\Theta_r(z)=\f{z^p}{r^{ps}}+\|a\|_\iy \,\f{z^q}{r^{qt}}\,\,\text{ and }\,\,\theta_r(z)=\f{z^{p-1}}{r^{ps}}+\|a\|_\iy \,\f{z^{q-1}}{r^{qt}}.
\end{split}\end{equation}

\,

{\bf [Proof of Theorem 1.4] }
For simplicity, set $x_0=0$. So in what follows we write $B^0_r\equiv B_r$ for $r>0$. Take any $r\in(0,R/2)$ and put 
$$u_b=u+b$$ for $b>0$. We consider a radial function $\zt\in C^\iy_c(B_{3r/2})$ with values in $[0,1]$ such that 
$$\zt|_{B_r}\equiv 1, \,\,\zt|_{\BR^n\s B_{2r}}\equiv 0\,\,\text{ and $\,\,|\n\zt|\ls 1/r\,$ in $\BR^n$.}$$ As a testing function, we use the function 
$$\vp(x)=\f{\zt^q(x)}{\theta_r(u_b(x))}.$$
Then we have that
\begin{equation}\begin{split}
0&\le\iint_{B_{2r}\times B_{2r}}\biggl[\,H_p(u_b(x)-u_b(y))(\vp(x)-\vp(y))K_{ps}(x,y)   \\
&\qquad\qquad\qquad+\fa(x,y)H_q(u_b(x)-u_b(y))(\vp(x)-\vp(y))K_{qt}(x,y)\biggr]\,dx\,dy \\
&\quad\qquad+2\int_{\BR^n\s B_{2r}}\int_{B_{2r}}
\biggl[\,H_p(u_b(x)-u_b(y))\vp(x)K_{ps}(x,y)   \\
&\qquad\qquad\qquad\qquad+\fa(x,y)H_q(u_b(x)-u_b(y))\vp(x)K_{qt}(x,y)\biggr]\,dx\,dy \\
&:=\cJ_1+\cJ_2.
\end{split}\end{equation}
Without loss of generality, we assume that 
$$u_b(x)\ge u_b(y)$$ for the estimate $\cJ_1$; for, by symmetry the other case $u_b(x)<u_b(y)$ can be done in the exactly same way. Then we have two possible cases (a) $u_b(x)\le 2 u_b(y)$ and (b) $u_b(x)> 2 u_b(y)$.

[ Case (a) : $u_b(y)\le u_b(x)\le 2 u_b(y)$ ]  From the mean value theorem, we note that
\begin{equation}\begin{split}
\zt(x)\ge\zt(y)\,\,&\Rightarrow\,\,\zt^q(x)-\zt^q(y)=q\int_{\zt(y)}^{\zt(x)}\tau^{q-1}d\tau\le q\zt^{q-1}(x)(\zt(x)-\zt(y)), \\
\zt(x)<\zt(y)\,\,&\Rightarrow\,\,\zt^q(x)-\zt^q(y)=q\int_{\zt(x)}^{\zt(y)}(-\tau^{q-1})d\tau\le q\zt^{q-1}(x)(\zt(x)-\zt(y)).
\end{split}\end{equation}
Since we see that
\begin{equation}\f{p-1}{\Theta_r(z)}\le\f{\theta'_r(z)}{\theta^2_r(z)}\le\f{q-1}{\Theta_r(z)},
\end{equation}
it follows from (6.7) and (6.8) that
\begin{equation}\begin{split}
&\vp(x)-\vp(y)=\f{\zt^q(x)-\zt^q(y)}{\theta_r(u_b(y)}+\zt^q(x)\biggl(\f{1}{\theta_r(u_b(x))}-\f{1}{\theta_r(u_b(y))}\biggr) \\
&\qquad\le\f{q\,\zt^{q-1}(x)(\zt(x)-\zt(y))}{\theta_r(u_b(y))}+\zt^q(x)\int_0^1\f{d}{d\tau}\biggl(\f{1}{\theta_r(\tau(u_b(x)-u_b(y))+u_b(y))}\biggr)\,d\tau \\
&\qquad\le\f{q\,\zt^{q-1}(x)(\zt(x)-\zt(y))}{\theta_r(u_b(y))}-(p-1)\f{\zt^q(x)(u_b(x)-u_b(y))}{\Theta_r(u_b(x))} \\
&\qquad\le\f{q\,\zt^{q-1}(x)\,|\zt(x)-\zt(y)|\,\,u_b(y)}{\Theta_r(u_b(y))}-\f{p-1}{2^q}\,\f{\zt^q(x)(u_b(x)-u_b(y))}{\Theta_r(u_b(x))} 
\end{split}\end{equation}
Using Young's inequality with indices $p'=\f{p}{p-1},p,\vep$, observing the fact that
\begin{equation*}(q-1)p'-q=\f{q-p}{p-1}\ge 0
\end{equation*} and applying the argument as in (5.12),
by (6.5), (6.6) and (6.9) we have that
\begin{equation}\begin{split}
\cJ_1&\le\iint_{B_{2r}\times B_{2r}}\biggl(\,\ld q\,\f{\theta_*(x,y;u_b(x)-u_b(y))\,\zt^{q-1}(x)\,|\zt(x)-\zt(y)|\,u_b(y)}{\Theta_r(u_b(y))} \\
&\qquad\qquad\qquad\qquad\qquad-\f{p-1}{\Ld 2^q}\,\f{\zt^q(x)\,\Theta_*(x,y;u_b(x)-u_b(y))}{\Theta_r(u_b(y))}\,\biggr)\,\f{dx\,dy}{|x-y|^{n}} \\
&\le\ld q\iint_{B_{2r}\times B_{2r}}\biggl(\,\f{\vep(u_b(x)-u_b(y))^p\zt^{(q-1)p'}(x)+c_\vep|\zt(x)-\zt(y)|^p u^p_b(y)}{|x-y|^{ps}\,\Theta_r(u_b(y))}   \\
&\qquad\qquad+a(x,y)\f{\vep(u_b(x)-u_b(y))^q\zt^q(x)+c_\vep|\zt(x)-\zt(y)|^q u^q_b(y)}{|x-y|^{qt}\,\Theta_r(u_b(y)) } \biggr)\,\f{dx\,dy}{|x-y|^{n}} \\
&\qquad\qquad\qquad-\f{c(p-1)}{\Ld 2^q}\iint_{B_{2r}\times B_{2r}}\f{\zt^q(x)\,\Theta(x,y;u_b(x)-u_b(y))}{\Theta_r(u_b(y))}\,\f{dx\,dy}{|x-y|^{n}} \\
&\le\ld q\,\vep\iint_{B_{2r}\times B_{2r}}\f{\zt^q(x)\,\Theta(x,y;u_b(x)-u_b(y))}{\Theta_r(u_b(y))}\,\f{dx\,dy}{|x-y|^{n}} \\
&\qquad\qquad+\ld q\, c_\vep\iint_{B_{2r}\times B_{2r}}\f{\Theta(x,y;(\zt(x)-\zt(y))u_b(y))}{\Theta_r(u_b(y))}\,\f{dx\,dy}{|x-y|^{n}} \\
&\qquad\qquad\qquad-\f{c(p-1)}{\Ld 2^q}\iint_{B_{2r}\times B_{2r}}\f{\zt^q(x)\,\Theta(x,y;u_b(x)-u_b(y))}{\Theta_r(u_b(y))}\,\f{dx\,dy}{|x-y|^{n}}.
\end{split}\end{equation}
Taking $\vep$ in the above (6.10) as
$$\vep=\f{c(p-1)}{\ld\Ld q 2^{q+1}},$$
we obtain that 
\begin{equation}\begin{split}
&\cJ_1\le\ld q\,c_\vep\iint_{B_{2r}\times B_{2r}}\f{\Theta(x,y;(\zt(x)-\zt(y))u_b(y))}{\Theta_r(u_b(y))}\,\f{dx\,dy}{|x-y|^{n}} \\
&\qquad-\f{c(p-1)}{\Ld 2^{q+1}}\iint_{B_{2r}\times B_{2r}}\f{\zt^q(x)\,\Theta(x,y;u_b(x)-u_b(y))}{\Theta_r(u_b(y))}\,\f{dx\,dy}{|x-y|^{n}} \\
&\le\ld q\,c_\vep\iint_{B_{2r}\times B_{2r}}\biggl(\f{r^{ps}(\zt(x)-\zt(y))^p}{|x-y|^{n+ps}}+\f{r^{qt}(\zt(x)-\zt(y))^q}{|x-y|^{n+qt}}\biggr)\,dx\,dy \\
&-\f{c(p-1)}{\Ld 2^{q+1}}\,r^{t(p-1)}\iint_{B_{2r}\times B_{2r}}\f{\zt^q(x)\,\Theta(x,y;u_b(x)-u_b(y))\,|x-y|^{t(p-1)}}{\Theta_r(u_b(y))\,r^{t(p-1)}}\,\f{dx\,dy}{|x-y|^{n+t(p-1)}}.
\end{split}\end{equation}
We now consider the monotone increasing function $\vartheta_1$ given by
$$\vartheta_1(z)=\f{z}{|x-y|^{(s-t)p}}+\|\fa\|_\iy \,z^{q/p},\,\,z\ge 0$$
and note that
\begin{equation}\begin{split}
\bigl(\ln u_b(x)-\ln u_b(y)\bigr)^p&=\biggl(\int_0^1\f{u_b(x)-u_b(y)}{\tau(u_b(x)-u_b(y))+u_b(y)}\,d\tau\biggr)^{p} \\
&\le\f{\ds\f{(u_b(x)-u_b(y))^p}{|x-y|^{pt}}}{\ds\f{u^p_b(y)}{r^{pt}}}\,\f{|x-y|^{pt}}{r^{pt}}.
\end{split}\end{equation}
Applying Lemma 3.1 to (6.12), we easily have that
\begin{equation}\begin{split}
\bigl(\ln u_b(x)-u_b(y)\bigr)^p&\le\left[\,\f{\vartheta_1\ds\biggl(\f{(u_b(x)-u_b(y))^p}{|x-y|^{pt}}\biggr)}{\vartheta_1\ds\biggl(\f{u^p_b(y)}{r^{pt}}\biggr)}+1\right]\f{|x-y|^{pt}}{r^{pt}} \\
&\le \f{\Theta(x,y;u_b(x)-u_b(y))\,|x-y|^{pt}}{\Theta_r(u_b(y))\,r^{pt}}+\f{|x-y|^{pt}}{r^{pt}}\\
&\ls \f{\Theta(x,y;u_b(x)-u_b(y))\,|x-y|^{t(p-1)}}{\Theta_r(u_b(y))\,r^{t(p-1)}}+\f{|x-y|^{pt}}{r^{pt}},
\end{split}\end{equation} because $|x-y|\le 2r$.
Putting (6.13) into (6.11), we have that
\begin{equation}\begin{split}
\cJ_1&\ls-r^{t(p-1)}\iint_{B_{2r}\times B_{2r}}\zt^q(x)\,\biggl|\ln\biggl(\f{u_b(x)}{u_b(y)}\biggr)\biggr|^p\,\f{dx\,dy}{|x-y|^{n+t(p-1)}} \\
&\qquad+\iint_{B_{2r}\times B_{2r}}\biggl(\f{|x-y|^{(1-s)p}}{r^{(1-s)p}}+\f{|x-y|^{(1-t)q}}{r^{(1-t)q}}+\f{|x-y|^t}{r^t}\biggr)\,\f{dx\,dy}{|x-y|^n} \\
&\ls-r^{t(p-1)}\iint_{B_{r}\times B_{r}}\,\biggl|\ln\biggl(\f{u_b(x)}{u_b(y)}\biggr)\biggr|^p\,\f{dx\,dy}{|x-y|^{n+t(p-1)}}+r^{n}.
\end{split}\end{equation}

[ Case (b) : $u_b(x)>2 u_b(y)$ ]  
It follows from the inequality with $\vep=\f{1}{2}(2^{p-1}-1)$ in Lemma 3.5 that
\begin{equation}\begin{split}
\vp(x)-\vp(y)&=\f{\zt^q(x)-\zt^q(y)}{\theta_r(u_b(x))}+\zt^q(y)\biggl(\f{1}{\theta_r(u_b(x))}-\f{1}{\theta_r(u_b(y))}\biggr) \\
&\le\f{\zt^q(x)-\zt^q(y)}{\theta_r(u_b(x))}+\zt^q(y)\biggl(\f{1}{\theta_r(2u_b(y))}-\f{1}{\theta_r(u_b(y))}\biggr) \\
&\le\f{\vep\zt^q(y)+c_\vep(\zt(x)-\zt(y))^q}{\theta_r(u_b(x))}-(1-2^{-p+1})\f{\zt^q(y)}{\theta_r(u_b(y))} \\
&\le c\,\f{(\zt(x)-\zt(y))^q}{\theta_r(u_b(x))}-\bigl(\f{1}{2}-2^{-p}\bigr)\f{\zt^q(y)}{\theta_r(u_b(y))}.
\end{split}\end{equation}
Applying the argument as in (5.12), by (6.15) we have that
\begin{equation}\begin{split}
\cJ_1&\le c\iint_{B_{2r}\times B_{2r}}\f{\theta_*(x,y;u_b(x)-u_b(y))\,(\zt(x)-\zt(y))^q}{\theta_r(u_b(x))}\,\f{dx\,dy}{|x-y|^{n}} \\
&\qquad\qquad-\bigl(\f{1}{2}-2^{-p}\bigr)\iint_{B_{2r}\times B_{2r}}\f{\zt^q(y)\,\theta_*(x,y;u_b(x)-u_b(y))}{\theta_r(u_b(y))}\,\f{dx\,dy}{|x-y|^{n}} \\
&\ls \iint_{B_{2r}\times B_{2r}}\f{\theta(x,y;u_b(x)-u_b(y))\,(\zt(x)-\zt(y))^q}{\theta_r(u_b(x))}\,\f{dx\,dy}{|x-y|^{n}} \\
&\quad-r^{t(p-1)}\iint_{B_{2r}\times B_{2r}}\f{\zt^q(y)\,\theta(x,y;u_b(x)-u_b(y))\,|x-y|^{t(p-1)}}{\theta_r(u_b(y))\,r^{t(p-1)}}\,\f{dx\,dy}{|x-y|^{n+t(p-1)}}.
\end{split}\end{equation}
It is easy to check that
\begin{equation}\f{\theta(x,y;u_b(x)-u_b(y))\,(\zt(x)-\zt(y))^q}{\theta_r(u_b(x))}
\le c\,\biggl(\,\f{r^{ps}}{|x-y|^{ps}}+\f{r^{qt}}{|x-y|^{qt}}\biggr)\,\f{|x-y|^q}{r^{q}}.
\end{equation}
Since $(\ln t)^p\le c\,(t-1)^{p-1}$ for $t>2$, we have that
\begin{equation}\begin{split}
\bigl(\ln u_b(x)-\ln u_b(y)\bigr)^p&\le c\,\biggl(\f{u_b(x)-u_b(y)}{u_b(y)}\biggr)^{p-1} \\
&=c\,\f{\ds\biggl(\f{u_b(x)-u_b(y)}{|x-y|^t}\biggr)^{p-1}}{\ds\biggl(\f{u_b(y)}{r^t}\biggr)^{p-1}}\,\f{|x-y|^{t(p-1)}}{r^{t(p-1)}}.
\end{split}\end{equation}
If we consider the increasing funtion
$$\vartheta_2(z)=\f{z}{|x-y|^{(s-t)p+t}}+\|\fa\|_\iy\,\f{z^{\f{q-1}{p-1}}}{|x-y|^{t}},\,\,z\ge 0,$$
then by Lemma 3.1 and (6.18) we have that
\begin{equation}\begin{split}
&\bigl|\ln u_b(x)-\ln u_b(y)\bigr|^p\le c\,\f{\vartheta_2\biggl(\,\ds\biggl(\f{u_b(x)-u_b(y)}{|x-y|^t}\biggr)^{p-1}\,\biggr)}{\vartheta_2\biggl(\,\ds\biggl(\f{u_b(y)}{r^t}\biggr)^{p-1}\,\biggr)}\,\f{|x-y|^{t(p-1)}}{r^{t(p-1)}} \\
&\quad\qquad\qquad\le c\,\,\f{\theta(x,y;u_b(x)-u_b(y))\,|x-y|^{t(p-1)}}{\theta_r(u_b(y))\,r^{t(p-1)}},
\end{split}\end{equation}
because $|x-y|\le 2r$. Thus it follows from (6.16), (6.17) and (6.19) that
\begin{equation}\begin{split}
\cJ_1&\ls-r^{t(p-1)}\iint_{B_{2r}\times B_{2r}}\zt^q(y)\,\biggl|\ln\biggl(\f{u_b(x)}{u_b(y)}\biggr)\biggr|^p\,\f{dx\,dy}{|x-y|^{n+t(p-1)}} \\
&\qquad+\iint_{B_{2r}\times B_{2r}}\biggl(\f{|x-y|^{q-ps}}{r^{q-ps}}+\f{|x-y|^{(1-t)q}}{r^{(1-t)q}}\biggr)\,\f{dx\,dy}{|x-y|^n} \\
&\le-r^{t(p-1)}\iint_{B_{r}\times B_{r}}\biggl|\ln\biggl(\f{u_b(x)}{u_b(y)}\biggr)\biggr|^p\,\f{dx\,dy}{|x-y|^{n+t(p-1)}}+r^{n}.
\end{split}\end{equation}

For the estimate of $\cJ_2$, we note that 
$u_b(x)-u_b(y)\le u_b(x)$ for $(x,y)\in B_{3r/2}\times(B_R\s B_{2r})$ and 
$$(u_b(x)-u_b(y))_+\le u_b(x)+u^-(y)$$ for $(x,y)\in B_{3r/2}\times(\BR^n\s B_R)$. This implies that
\begin{equation}\begin{split}
\cJ_2&\le 2\int_{B_{3r/2}}\int_{\BR^n\s B_{2r}}\f{\theta(x,y;u_b(x))}{\theta_r(u_b(x))}\,\f{dy\,dx}{|x-y|^{n}} \\
&\qquad\qquad\qquad+2\int_{B_{3r/2}}\int_{\BR^n\s B_R}\f{\theta(x,y;u^-(y))}{\theta_r(u_b(x))}\,\f{dy\,dx}{|x-y|^{n}} \\
&\ls\int_{B_{3r/2}}\int_{\BR^n\s B_{2r}}\f{r^{qt}}{|x-y|^{n+qt}}\,dy\,dx \\
&\qquad+\f{r^{ps}}{b^{p-1}}\int_{B_{3r/2}}\int_{\BR^n\s B_R}\biggl(\,\f{[u^-(y)]^{p-1}}{|y-x_0|^{n+ps}}+\|a\|_\iy\,\f{[u^-(y)]^{q-1}}{|y-x_0|^{n+qt}}\biggr)\,dy\,dx \\
&\ls r^{n}+\f{r^{n}}{b^{p-1}}\biggl(\f{r}{R}\biggr)^{ps}\,\cT_R(u^-;x_0),
\end{split}\end{equation}
because $u(x)\ge 0$ for $x\in B_R$ and 
$|x-y|\ge|y-x_0|/4$ for $(x,y)\in B_{3r/2}\times(\BR^n\s B_{2r})$.
Since $x,y\in B_r$, it thus follows from (6.6), (6.14), (6.20) and (6.21) that
\begin{equation*}\begin{split}
\iint_{B_r\times B_r}\biggl|\ln\biggl(\f{u_b(x)}{u_b(y)}\biggr)\biggr|^p\,\f{dx\,dy}{|x-y|^{n+t(p-1)}}
&\ls r^{n-t(p-1)}\biggl[1+\f{1}{b^{p-1}}\biggl(\f{r}{R}\biggr)^{ps}\,\cT_R(u^-;x_0)\biggr] \\
\end{split}\end{equation*} 
Hence we complete the proof. \qed
 
\begin{lemma} Suppose that $qt\le ps\le n\,$ for $0<s,t<1$ and $1<p\le q\le p_s^*$, and let $g\in L^{q-1}_{qt}(\BR^n)$. 
If $\,w=(v\vee 0)\wedge d\,$ for any $a,b,d>0$ where
$$v=\ln\,\biggl(\f{a+b}{u(x)+b}\biggr)$$
and $u\in\fD^-_g(\Om)$ is a weak supersolution of the nonlocal double phase equation $(1.5)$ with $u\ge 0$ in $B^0_R\subset\Om$, then we have that $w\in\bmo(B^0_R)$ and moreover
$$\|w\|_{\bmo(B^0_R)}\ls\biggl(1+\f{1}{b^{p-1}}\biggl(\f{r}{R}\biggr)^{ps}\,\cT_R(u^-;x_0)\biggr)^{\f{1}{p}}<\iy.$$
\end{lemma}

\pf We note that $\bigl||v(x)|-|v|_{B_r}\bigr|\le|v(x)-v_{B_r}|$ and 
$$W^{s,p}(\Om)\subset W^{t_{p},p}(\Om)$$
because $0<t_{p}:=t(1-1/p)<t<s$. Since $pt_p<ps\le n$, by the fractional Poincar\'e inequality (6.1), (6.3), (6.4)
and Theorem 1.4 we have that
\begin{equation*}\begin{split}
\biggl(\intavg_{B_r}\bigl||v|-|v|_{B_r}\bigr|\,dx\biggr)^p&\le\biggl(\intavg_{B_r}|v-v_{B_r}|\,dx\biggr)^p
\le \intavg_{B_r}|v-v_{B_r}|^p\,dx\\
&\le\f{1}{|B_r|}\,|B_r|^{\f{p\,t_{p}}{n}}[v]^p_{W^{t_{p},p}(B_r)}\\
&\ls\f{1}{|B_r|}\,|B_r|^{\f{p\,t_{p}}{n}}\,r^{n-p\,t_{p}}\biggl[1+\f{1}{b^{p-1}}\biggl(\f{r}{R}\biggr)^{ps}\,\cT_R(u^-;x_0)\biggr] \\
&\ls 1+\f{1}{b^{p-1}}\biggl(\f{r}{R}\biggr)^{ps}\,\cT_R(u^-;x_0)<\iy.
\end{split}\end{equation*}
Hence we are done. \qed

\,

Next we furnish a fundamental inequality related with the Sobolev-Poincar\'e inequality (Lemma 2.3), which is useful in obtaining the interior H\"older regularity of weak solutions to nonlocal double phase equation (1.5). For $x,y\in\BR^n$, $r>0$ and $z_1,z_2\ge 0$, we consider two functions
$$\Theta_r(z_1,z_2)=\f{z_1^p}{r^{ps}}+\|\fa\|_\iy\,\f{z_2^q}{r^{qt}}\,\,
\text{ and }\,\,\Theta(x,y;z_1,z_2)=\f{z_1^p}{|x-y|^{ps}}+\|\fa\|_\iy\,\f{z_2^q}{|x-y|^{qt}}.$$

\begin{lemma} Let $qt\le ps\le n$ for $0<s,t<1<p\le q\le p_s^*$, and $u,v\in W^{s,p}(B^0_r)\cap W^{t,q}(B^0_r)$ for $r>0$.

$(a)$ If $ps<n$, then there is a constant $c=c(n,p,s)>0$ such that
\begin{equation*}\begin{split}
\biggl(\,\intavg_{B^0_r}[\Theta_r(|u(x)|,|v(x)|)]^\bt\,dx\biggr)^{\f{1}{\bt}}&\le c\intavg_{B^0_r}\int_{B^0_r}\f{\Theta(x,y;|u(x)-u(y)|,|v(x)-v(y)|)}{|x-y|^n}\,dx\,dy \\
&\qquad\qquad\qquad\qquad+c\intavg_{B^0_r}\Theta_r(|u(x)|,|v(x)|)\,dx
\end{split}\end{equation*} 
where $\bt=q_t^*/q>1$.

$(b)$  If $ps=n$, then there is a constant $c=c(n,p,s)>0$ such that
\begin{equation*}\begin{split}
\biggl(\,\intavg_{B^0_r}[\Theta_r(|u(x)|,|v(x)|)]^{\bt_1}\,dx\biggr)^{\f{1}{\bt_1}}&\le c\intavg_{B^0_r}\int_{B^0_r}\f{\Theta(x,y;|u(x)-u(y)|,|v(x)-v(y)|)}{|x-y|^n}\,dx\,dy \\
&\qquad\qquad\qquad\qquad+c\intavg_{B^0_r}\Theta_r(|u(x)|,|v(x)|)\,dx
\end{split}\end{equation*} 
where $\bt_1=q_{t_1}^*/q>1$ for $t_1\in(0,1)$ so that $qt_1\le ps_1<n$ for $s_1\in(0,1)$.
\end{lemma}

\pf (a) Since $ps\ge qt$, we note that
$$1<\f{q_t^*}{q}\le\f{p_s^*}{p}.$$
Since the function 
$$\fA_f(p):=\biggl(\,\intavg_E |f(x)|^p\,dx\biggr)^{\f{1}{p}}$$ is increasing on $(0,\iy)$ for $f\in\cF(\BR^n)$
and a measurable set $E\subset\BR^n$ with $0<|E|<\iy$, it follows from Lemma 4.2 and Jensen's inequality that
\begin{equation}\begin{split}
\biggl(\,\intavg_{B^0_r}[\Theta_r(|u(x)|,&|v(x)|)]^{\bt}\,dx\biggr)^{\f{1}{\bt}} 
\le c\,\biggl(\,\intavg_{B^0_r}[\Theta_r(|u(x)-u_{B^0_r}|,|v(x)-v_{B^0_r}|)]^{\bt}\,dx\biggr)^{\f{1}{\bt}} \\
&\qquad\qquad\qquad\qquad+c\,\Theta_r(|u_{B^0_r}|,|v_{B^0_r}|)  \\
&\le c\,\biggl[\biggl(\,\intavg_{B^0_r}\biggl|\f{u-u_{B^0_r}}{r^s}\biggr|^{p_s^*}\,dx\biggr)^{\f{p}{p_s^*}}
+\|\fa\|_\iy\biggl(\,\intavg_{B^0_r}\biggl|\f{v-v_{B^0_r}}{r^t}\biggr|^{q_t^*}\,dx\biggr)^{\f{q}{q_t^*}}\biggr] \\
&\qquad\qquad\qquad\qquad+c\,\Theta_r(|u_{B^0_r}|,|v_{B^0_r}|) \\
&\le c\intavg_{B^0_r}\int_{B^0_r}\f{\Theta(x,y;|u(x)-u(y)|,|v(x)-v(y)|)}{|x-y|^n}\,dx\,dy \\
&\qquad\qquad\qquad\qquad+c\intavg_{B^0_r}\Theta_r(|u(x)|,|v(x)|)\,dx.
\end{split}\end{equation}

(b) It can be done in a similar way. 
Hence we complete the proof. \qed

\,

\section{ The Supremum and infimum estimate of weak supersolutions }

In this section, we obtain a upper bound for the supremum estimate of the local ratio of the lowerer level set of weak supersolutions of the nonlocal double phase equation (1.5) with $u\ge 0$ in $B^0_R\subset\Om$. Also we get a lower bound of the infimum estimate of its weak supersolutions.

\begin{lemma} Let $qt\le ps\le n$ for $0<s,t<1$ and $1<p\le q\le p_s^*$, $g\in L^{q-1}_{qt}(\BR^n)$ and let $a\ge 0$. 
Suppose that $u\in\fD^+_g(\Om)$ is a weak supersolution of the nonlocal double phase equation $(1.5)$ with $u\ge 0$ in $B^0_R\subset\Om$. If there is some $\nu\in(0,1)$ such that
\begin{equation}\inf_{r\in(0, R/2)}\f{|B^0_r\cap\{u\ge a\}}{|B^0_r|}\ge\nu,
\end{equation}
then for each $r\in(0,R/2)$ there is a constant $c_0=c(n,s,t,p,q,\ld,\Ld,\Om,\fa)>0$ such that
$$\sup_{\dt\in(0,1/4)}\ln\bigl(\ds\f{1}{2\dt}\bigr)\f{|B^0_r\cap\{u\le 2\dt a-b\}|}{|B^0_r|}\le\f{c_0}{\nu},$$
where $$b=\f{1}{2}\biggl(\f{r}{R}\biggr)^{\f{ps}{p-1}}\,\cT^{\f{1}{p-1}}_R(u^-;x_0).$$
\end{lemma}

\pf Take any $r\in(0,R/2)$. For $\dt\in(0,1/4)$, we define the function $w$ by
$$w=\biggl(\ln\,\biggl(\f{a+b}{u+b}\biggr)\vee 0\biggr)\wedge\ln\biggl(\f{1}{2\dt}\biggr).$$
Then we easily see from Lemma 6.1 that
\begin{equation}\intavg_{B^0_r}|w(x)-w_{B^0_r}|\,dx\le c_0.
\end{equation}
From the definition of $w$, we know that
$$B^0_r\cap\{w=0\}=B^0_r\cap\{u\ge a\}.$$
Then by (7.1) we have that
\begin{equation}|B^0_r\cap\{w=0\}|\ge\nu|B^0_r|.
\end{equation}
Then it follows from (7.3) that
\begin{equation}\begin{split}
\ln\,\biggl(\f{1}{\dt}\biggr)&=\f{1}{|B^0_r\cap\{w=0\}|}\int_{B^0_r\cap\{w=0\}}\biggl[\ln\,\biggl(\f{1}{\dt}\biggr)-w(x)\biggr]\,dx \\
&\le\f{1}{\nu}\,\biggl[\ln\,\biggl(\f{1}{\dt}\biggr)-w_{B^0_r}\biggr].
\end{split}\end{equation}
Applying (7.2) and integrating (7.4) on $B^0_r\cap\{w=\ln(1/(2\dt))\}$, we obtain that
\begin{equation*}\begin{split}
\biggl|B^0_r\cap\biggl\{w=\ln\,\biggl(\f{1}{2\dt}\biggr)\biggr\}\biggr|\ln\,\biggl(\f{1}{2\dt}\biggr)
&\le\f{1}{\nu}\int_{B^0_r}|w(x)-w_{B^0_r}|\,dx \le\f{c_0}{\nu}\,|B^0_r|.
\end{split}\end{equation*}
This makes it possible to get the inequality
$$\ln\bigl(\ds\f{1}{2\dt}\bigr)\,|B^0_r\cap\{u+b\le 2\dt a\}|\le\f{c_0}{\nu}|B^0_r|.$$
Hence we complete the proof. \qed

\begin{lemma} Let $qt\le ps\le n$ for $0<s,t<1$ and $1<p\le q\le p_s^*$, $g\in L^{q-1}_{qt}(\BR^n)$ and let $a\ge 0$. 
Suppose that $u\in\fD^+_g(\Om)$ is a weak supersolution of the nonlocal double phase equation $(1.5)$ with $u\ge 0$ in $B^0_R\subset\Om$. If there is some $\nu\in(0,1)$ such that
\begin{equation*}\inf_{r\in(0, R/2)}\f{|B^0_r\cap\{u\ge a\}}{|B^0_r|}\ge\nu,
\end{equation*}
then there is a constant $\dt=\dt(n,s,t,p,q,\ld,\Ld,\Om,\fa)\in(0,1/4)$ such that
\begin{equation*}\inf_{r\in(0,R/2)}\inf_{B_r(x)} u\ge \dt a-2 b,
\end{equation*}
where $$b=\f{1}{2}\biggl(\f{r}{R}\biggr)^{\f{ps}{p-1}}\,\cT^{\f{1}{p-1}}_R(u^-;x_0).$$
\end{lemma}

\pf First of all, we consider the case $qt\le ps<n$ and $1<p\le q\le p_s^*$.
In what follows, we write $B_r:=B^0_r$ for $r>0$. Since $u\ge 0$ in $B_R$, without loss of generality we may assume that 
\begin{equation} \dt a\ge 2 b. \end{equation}
Set $w=(h-u)_+$ for $h\in(\dt a,2\dt a)$ where $\dt>0$ is a constant to be selected later.
Take any $r\in(0,R/2)$. For $\vr\in(r,2r)$, we choose a testing function $$\vp(x)=w(x)\zt^q(x)$$ where $\zt\in C_c^\iy(B_\vr)$ is a function satisfying that $0\le\zt\le 1$ and $|\n\zt|\ls 1/\vr$ in $\BR^n$.
Then we have that
\begin{equation}\begin{split}
0&\le\iint_{B_\vr\times B_\vr}\biggl(H_p(u(x)-u(y))(\vp(x)-\vp(y)) K_{ps}(x,y) \\
&\qquad\qquad+\fa(x,y)H_q(u(x)-u(y))(\vp(x)-\vp(y)) K_{qt}(x,y)\biggr)\,dx\,dy \\
&\quad+2\iint_{(\BR^n\s B_\vr)\times B_\vr}\biggl(H_p(u(x)-u(y))\,\vp(x) K_{ps}(x,y) \\
&\quad\qquad\qquad\qquad+\fa(x,y)H_q(u(x)-u(y))\,\vp(x) K_{qt}(x,y)\biggr)\,dx\,dy \\
&:=\cI_1(u,\vp)+\cI_2(u,\vp).
\end{split}\end{equation}
Splitting $\cI_2(u,\vp)/2$ into two parts leads us to get that
\begin{equation*}\begin{split}
\f{1}{2}\,\cI_2(u,\vp)&=\int_{(\BR^n\s B_\vr)\cap\{y:u(y)<0\}}\int_{B_\vr}\biggl(H_p(u(x)-u(y))\,\vp(x) K_{ps}(x,y) \\
&\quad\qquad\qquad\qquad+\fa(x,y)H_q(u(x)-u(y))\,\vp(x) K_{qt}(x,y)\biggr)\,dx\,dy \\
&+\int_{(\BR^n\s B_\vr)\cap \{y:u(y)\ge 0\}}\int_{B_\vr}\biggl(H_p(u(x)-u(y))\,\vp(x) K_{ps}(x,y) \\
&\quad\qquad\qquad\qquad+\fa(x,y)H_q(u(x)-u(y))\,\vp(x) K_{qt}(x,y)\biggr)\,dx\,dy \\
&:=\cI_{2,1}(u,\vp)+\cI_{2,2}(u,\vp).
\end{split}\end{equation*}
Then we have that
\begin{equation*}\begin{split}
\cI_{2,1}(u,\vp)&\le h|B_\vr\cap\{x:u(x)<h\}| \\
&\sup_{x\in\supp(\zt)}\int_{\BR^n\s B_\vr}\biggl(\,[h+u^-(y)]^{p-1}K_{ps}(x,y)+\|\fa\|_\iy[h+u^-(y)]^{q-1}K_{qt}(x,y)\biggr)\,dy 
\end{split}\end{equation*} and
\begin{equation*}\begin{split}
\cI_{2,2}(u,\vp)&\le h^p|B_{\vr}\cap\{x:u(x)<h\}|\sup_{x\in\supp(\zt)}\int_{\BR^n\s B_\vr}K_{ps}(x,y)\,dy \\
&\qquad+h^q|B_{\vr}\cap\{x:u(x)<h\}|\,\,\|\fa\|_\iy\sup_{x\in\supp(\zt)}\int_{\BR^n\s B_\vr}K_{qt}(x,y)\,dy.
\end{split}\end{equation*}
Thus this leads us to obtain that
\begin{equation}\begin{split}
&\cI_2(u,\vp)\le 4\,h|B_\vr\cap\{x:u(x)<h\}| \\
&\qquad\sup_{x\in\supp(\zt)}\int_{\BR^n\s B_\vr}\biggl(\,[h+u^-(y)]^{p-1}K_{ps}(x,y)+\|\fa\|_\iy[h+u^-(y)]^{q-1}K_{qt}(x,y)\biggr)\,dy. \\
\end{split}\end{equation}
Now we set $\xi=\zt^{q/p}$. As in the proof of the Caccioppoli type estimate in (5.8), we have that
\begin{equation*}\begin{split}
\cI_1(u,\vp)&\ls-\iint_{B_\vr\times B_\vr}\biggl(\,|w(x)\xi(x)-w(y)\xi(y)|^p K_{ps}(x,y) \\
&\qquad\qquad\qquad+\fa(x,y)\,|w(x)\zt(x)-w(y)\zt(y)|^q K_{qt}(x,y)\biggr)\,dx\,dy \\
&\quad+\iint_{B_\vr\times B_\vr}\biggl(\,(w(x)\vee w(y))^p|\xi(x)-\xi(y)|^p K_{ps}(x,y)  \\
&\qquad\qquad\qquad+\fa(x,y)\,(w(x)\vee w(y))^q|\zt(x)-\zt(y)|^q K_{qt}(x,y)\biggr)\,dx\,dy.
\end{split}\end{equation*}
Applying again the argument as in (5.12), we obtain that
\begin{equation}\begin{split}
\cI_1(u,\vp)&\ls-\iint_{B_\vr\times B_\vr}\biggl(\,|w(x)\xi(x)-w(y)\xi(y)|^p K_{ps}(x,y) \\
&\qquad\qquad\qquad+\|\fa\|_\iy\,|w(x)\zt(x)-w(y)\zt(y)|^q K_{qt}(x,y)\biggr)\,dx\,dy \\
&\quad+\iint_{B_\vr\times B_\vr}\biggl(\,(w(x)\vee w(y))^p|\xi(x)-\xi(y)|^p K_{ps}(x,y)  \\
&\qquad\qquad\qquad+\|\fa\|_\iy\,(w(x)\vee w(y))^q|\zt(x)-\zt(y)|^q K_{qt}(x,y)\biggr)\,dx\,dy.
\end{split}\end{equation}
Then it follows from (7.6), (7.7) and (7.8) that
\begin{equation}\begin{split}
&[w\xi]^p_{W^{s,p}(B_\vr)}+\|\fa\|_\iy\,[w\zt]^q_{W^{t,q}(B_\vr)}\le 8\,h\,|B_\vr\cap\{u<h\}|  \\
&\quad\times\sup_{x\in\supp(\zt)}\int_{\BR^n\s B_\vr}\biggl(\,[h+u^-(y)]^{p-1}K_{ps}(x,y)+\|\fa\|_\iy[h+u^-(y)]^{q-1}K_{qt}(x,y)\biggr)\,dy \\
&\qquad\qquad+\iint_{B_\vr\times B_\vr}\biggl(\,(w(x)\vee w(y))^p|\xi(x)-\xi(y)|^p K_{ps}(x,y)  \\
&\qquad\qquad\qquad\qquad+\|\fa\|_\iy\,(w(x)\vee w(y))^q|\zt(x)-\zt(y)|^q K_{qt}(x,y)\biggr)\,dx\,dy \\
&\quad:=h\,|B_\vr\cap\{u<h\}|\,\cA_\vr(u,\zt,h)+\cB_\vr(w,\zt).
\end{split}\end{equation}
We are now applying Lemma 3.3 to our estimate (7.9). For $k=0,1,2,\cdots,$ we set
\begin{equation}h_k=\dt a+2^{-(k+1)}\dt a,\,\,\vr_k=r+2^{-k}r,\,\,\overline\vr_k=\f{\vr_{k+1}+\vr_k}{2}
\end{equation}
Here we note that $r<\vr_k,\overline\vr_k<2r$ and $h_k-h_{k+1}\ge 2^{-(k+3)}h_k$, and 
$$h_0=\f{3}{2}\,\dt a\le 2\dt a-\f{1}{2}\biggl(\f{r}{R}\biggr)^{\f{ps}{p-1}}\cT^{\f{1}{p-1}}_R(u^-;x_0).$$
Then we we see that
$$\{u<h_0\}\subset\biggl\{u\le 2\dt a-\f{1}{2}\biggl(\f{r}{R}\biggr)^{\f{ps}{p-1}}\cT^{\f{1}{p-1}}_R(u^-;x_0)\biggr\}.$$
Thus by Lemma 7.1 we have that
\begin{equation}\f{|B_r\cap\{u<h_0\}|}{|B_r|}\le\f{c_0}{\nu\,\ln\,\bigl(\f{1}{2\dt}\bigr)}.
\end{equation}
For each $k=0,1,2,\cdots,$ we set $B_k=B_{\vr_k}$ and let $\zt_k\in C^\iy_c(B_{\overline\vr_k})$ be a function with $\zt_k\equiv 1$ on $B_{\vr_{k+1}}$ such that 
$$0\le\zt_k\le 1\,\,\text{ and }\,\,|\n\zt_k|\le c/\vr_k\,\,\text{ in $\BR^n$.}$$
From (7.9) and  the fact that
$$|y-x|\ge|y-x_0|-|x-x_0|\ge\biggl(1-\f{\overline\vr_k}{\vr_k}\biggr)\,|y-x_0|\ge 2^{-k-3}|y-x_0|$$
for any $y\in\BR^n\s B_k$ and $x\in B_{\overline\vr_k}$, by (7.5) and (7.10) we obtain the estimate
\begin{equation}\begin{split}
&\cA_{\vr_k}(u,\zt_k,h_k) \\
&=\sup_{x\in\supp(\zt_k)}\int_{\BR^n\s B_k}\biggl(\,[h_k+u^-(y)]^{p-1}K_{ps}(x,y)+\|\fa\|_\iy[h_k+u^-(y)]^{q-1}K_{qt}(x,y)\biggr)\,dy \\
&\ls 2^{k(n+ps)}\int_{\BR^n\s B_k}\biggl(\,\f{[h_k+u^-(y)]^{p-1}}{|y-x_0|^{n+ps}}+\|\fa\|_\iy\f{[h_k+u^-(y)]^{q-1}}{|y-x_0|^{n+qt}}\biggr)\,dy \\
&\ls 2^{k(n+ps)}\theta_r^*(h_k)+2^{k(n+ps)}\int_{\BR^n\s B_R}\biggl(\,\f{[u^-(y)]^{p-1}}{|y-x_0|^{n+ps}}+\|\fa\|_\iy\f{[u^-(y)]^{q-1}}{|y-x_0|^{n+qt}}\biggr)\,dy \\
&\ls 2^{k(n+ps)}\theta_r^*(h_k)+2^{k(n+ps)}\,r^{-ps}\,\biggl(\f{r}{R}\biggr)^{ps}\cT_R(u^-;x_0) \ls 2^{k(n+ps)}\theta_r^*(h_k),
\end{split}\end{equation}
because $ps\ge qt$.
Also, we have the following estimate
\begin{equation}\begin{split}
\cB_{\vr_k}(w_k,\zt_k)
&\ls\iint_{B_k\times B_k}\biggl(\,(w_k(x)\vee w_k(y))^p\f{|\xi_k(x)-\xi_k(y)|^p}{|x-y|^{n+ps}}  \\
&\qquad\qquad\qquad+\|\fa\|_\iy\,(w_k(x)\vee w_k(y))^q\f{|\zt_k(x)-\zt_k(y)|^q}{|x-y|^{n+qt}} \biggr)\,dx\,dy \\
&\ls h_k^p\int_{B_k}\int_{B_k\cap\{u<h_k\}}\biggl(\,\f{\sup_{\BR^n}|\n\xi_k|^p}{|x-y|^{n+ps-p}}
+\f{\|\fa\|_\iy\sup_{\BR^n}|\n\zt_k|^q}{|x-y|^{n+qt-q}}  \biggr)\,dx\,dy  \\
&\ls h_k^p\,\biggl(\f{2^k}{r}\biggr)^p\int_{B_k\cap\{u<h_k\}}\int_{B_k}\f{1}{|y-x_0|^{n+ps-p}}\,dy\,dx \\
&\qquad+h_k^q\,\biggl(\f{2^k}{r}\biggr)^q\int_{B_k\cap\{u<h_k\}}\int_{B_k}\f{1}{|y-x_0|^{n+qt-q}}\,dy\,dx \\
&\ls 2^{qk}\,\Theta_r(h_k)\,|B_k\cap\{u<h_k\}|.
\end{split}\end{equation}
Thus it follows from (7.9), (7.12) and (7.13) that
\begin{equation*}\begin{split}
&\iint_{B_k\times B_k}\f{\Theta(x,y;|w_k(x)\xi_k(x)-w_k(y)\xi_k(y)|,|w_k(x)\zt_k(x)-w_k(y)\zt_k(y)|)}{|x-y|^n}\,dx\,dy \\
&\qquad\qquad\qquad\qquad\qquad\qquad\ls 2^{k(n+ps+q)}\,\Theta_r(h_k)\,|B_k\cap\{u<h_k\}|.
\end{split}\end{equation*}
Also, we easily obtain that
\begin{equation*}\int_{B_k}\Theta_{\vr_k}(|w_k\xi_k|,|w_k\zt_k|)\,dx\ls \Theta_r(h_k)\,|B_k\cap\{u<h_k\}|.
\end{equation*}
We apply (a) of Lemma 6.2 with index $\bt=\f{n}{n-qt}$ to derive the following inequality
\begin{equation}\begin{split}
&\biggl(\,\intavg_{B_{k+1}}[\Theta_{\vr_k}(|w_k|)]^\bt\,dx\biggr)^{\f{1}{\bt}} \\
&\quad\ls\intavg_{B_k}\int_{B_k}\f{\Theta(x,y;|w_k(x)\xi_k(x)-w_k(y)\xi_k(y)|,|w_k(x)\zt_k(x)-w_k(y)\zt_k(y)|)}{|x-y|^n}\,dx\,dy \\
&\qquad\qquad\qquad\qquad\qquad+\intavg_{B_k}\Theta_{\vr_k}(|w_k\xi_k|,|w_k\zt_k|)\,dx \\
&\quad\ls 2^{k(n+ps+q)}\,\Theta_r(h_k)\,\f{|B_k\cap\{u<h_k\}|}{|B_k|}.
\end{split}\end{equation}
For each $k=0,1,2,\cdots,$ we observe that
\begin{equation}w_k:=(h_k-u)_+\ge(h_k-h_{k+1})\mathbbm{1}_{\{u<h_{k+1}\}}\ge 2^{-k-3} h_k\mathbbm{1}_{\{u<h_{k+1}\}}.
\end{equation}
For $k=0,1,2,\cdots,$ we set
$$N_k=\f{|B_k\cap\{u<h_k\}|}{|B_k|}.$$
By (7.14) and (7.15), we obtain that
\begin{equation*}\begin{split}
\f{2^{-(k+3)q}}{2^{ps}}\,\Theta_r(h_k)\,N_{k+1}^{\f{1}{\bt}}&\le\biggl(\,\intavg_{B_{k+1}}[\Theta_{\vr_k}(|w_k|)]^\bt\,dx\biggr)^{\f{1}{\bt}} \le c\,2^{k(n+ps+q)}\,\Theta_r(h_k)\,N_k,
\end{split}\end{equation*}
which gives that
$$N_{k+1}^{\f{1}{\bt}}\le c\,2^{k(n+ps+2q)}\,N_k.$$
This leads us to get that
$$N_{k+1}\le c^\bt\,2^{k\bt(n+ps+2q)}\,N_k^{1+\f{qt}{n-qt}}.$$
From (7.11), we see that 
$$N_0\le\f{c_0}{\nu\,\ln\,\bigl(\f{1}{2\dt}\bigr)}.$$
We now apply Lemma 3.3 with
$$d_0=c^\bt,\,\,e_0=2^{\bt(n+ps+2q)}>1\,\,\text{ and }\,\,\e=\f{qt}{n-qt}>0.$$
If we choose a small $\dt$ depending only on $n,s,t,p,q,\ld,\Ld,\Om$ and $\fa$ so that
$$0<\dt:=\f{1}{2}\exp\,\biggl(-\f{c_0 c^{\f{\bt(n-qt)}{qt}} 2^{\f{n(n-qt)}{q^2t^2}(n+ps+2q)}}{\nu}\biggr)<\f{1}{4},
$$ then we have that
$$N_0\le c^{-\f{\bt(n-qt)}{qt}}(2^{n+ps+2q})^{-\f{(n-qt)^2}{q^2 t^2}}.$$
Thus we conclude that $\lim_{k\to\iy}N_k=0$, which implies that
$$\inf_{B_r}u\ge\dt a.$$ 

Secondly, we consider the case $qt\le ps=n$ and $1<p\le q<\iy$. Using (b) of Lemma 6.2, we can obtain the required result in a similar way. Hence we complete the proof. \qed

\section{ Nonlocal weak Harnack inequalities }

In this section, first of all, we mention a variant of the Krylov-Safonov covering theorem \cite{KS,KSh} which is useful in obtaining nonlocal weak Harnack inequalities. If $E\subset\BR^n$ is a bounded set and $\cC_E=\{B_r(x)\subset\BR^n:x\in E\}$, then it follows from Vitali's covering lemma that there is a countable pairwise disjoint subcovering $\fC_E=\{B_{r_k}(x_k)\}_{k\in\BN}$ of $\cC_E$ such that
$$E\subset\bigcup_{k\in\BN}B_{5 r_k}(x_k).$$
For $\vr>0$, $\ap\in(0,1)$ and a measurable subset $E$ of a ball $B_r(x_0)$, we define the set
\begin{equation}E_\ap^\vr:=\bigcup_{\rho\in(0,\vr)}\bigl\{B_\rho(x)\cap B_r(x_0):|E\cap B_{5\rho}(x)|>\ap|B_\rho(x)|,\,x\in B_r(x_0)\bigr\}
\end{equation}
The following version of the Krylov-Safonov covering theorem no longer depends on the treshold radius $\vr$, whose proof can be obtained from that of \cite{KSh} ( refer to \cite{K1} for the detailed proof ).

\begin{lemma} If $\vr>0$, $\ap\in(0,1)$ and $E$ is a measurable subset of $B_r(x_0)$, then either $E^\vr_\ap=B_r(x_0)$ or
$$|E^\vr_\ap|\ge\f{1}{(1+\ep)^n\ap}\,|E|\,\,\text{ for any $\ep\in(0,1)$. }$$
\end{lemma}

Combining Lemma 8.1 with the estimates obtained in Section 7 lead us to derive nonlocal weak Harnack inequalities of weak supersolutions to the nonlocal diuble phase equation (1.4) as follows.

\,

{\bf [Proof of Theorem 1.6]} Take any $r\in(0,R/2)$. For simplicity, we write $B_r:=B^0_r$ for $r>0$. For $k\in\BN_0$ and $\tau>0$, we consider the set
$$U^k_\tau=\biggl\{y\in B_{2r}:u(y)\ge\tau\dt^k-\f{2 b}{1-\dt}\biggr\}$$
where $\dt\in(0,1)$ is the constant in Lemma 7.2 and $b>0$ is the constant given by
$$b=\f{1}{2}\biggl(\f{r}{R}\biggr)^{\f{ps}{p-1}}\,\cT^{\f{1}{p-1}}_R(u^-;x_0).$$
Then it is obvious that $U^{k-1}_\tau\subset U^k_\tau$ for all $k\in\BN$. Choose any point $x\in B_{2r}$ so that $B_\rho(x)\cap B_{2r}\subset E^\vr_\ap$ with $E=U^{k-1}_\tau$. From (8.1), we have that
$$|U^{k-1}_\tau\cap B_{5\rho}(x)|>\ap\,|B_{\rho}(x)|=\f{\ap}{5^n}\,|B_{5\rho}(x)|.$$
Applying Lemma 7.2 with $a=\tau\dt^{k-1}-2b/(1-\dt)$, we can derive that
$$u\ge\dt\,\biggl(\tau\dt^{k-1}-\f{2b}{1-\dt}\biggr)-2b=\tau\dt^k-\f{2b}{1-\dt}\,\,\text{ in $B_{5\rho}(x)\cap B_{2r}$,}$$
and thus we have that $E^\vr_\ap\subset U^k_\tau$. Hence Lemma 8.1 makes it possible to obtain that either $U^k_\tau=B_{2r}$ or
$$|U^k_\tau|\ge\f{1}{(1+\ep)^n\ap}\,|U^{k-1}_\tau|\,\,\text{ for any $k\in\BN$ and $\ep\in(0,1)$. }$$
We now choose some $\ep_0\in(0,1)$ such that
$$1/16<(1+\ep_0)^n\ap<1.$$
Then we claim that $U^N_\tau=B_{2r}$ for some $N\in\BN$ with
$$|U^0_\tau|>\bigl((1+\ep_0)^n\ap\bigr)^N|B_{2r}|;$$
indeed, if $U^N_\tau\neq B_{2r}$, then we have that
$$|U^N_\tau|\ge\f{1}{(1+\ep_0)^n\ap}\,|U^{N-1}_\tau|,$$
and so $U^k_\tau\neq B_{2r}$ for any $k=1,2,\cdots,N$. This implies that
$$|U^N_\tau|\ge\f{1}{(1+\ep_0)^n\ap}\,|U^{N-1}_\tau|\ge\cdots\ge\biggl(\f{1}{(1+\ep_0)^n\ap}\biggr)^N|U^0_\tau|>|B_{2r}|,$$
which gives a contradiction. Thus it follows from the calim $U^N_\tau=B_{2r}$  that
\begin{equation}u\ge\tau\dt^N-\f{2b}{1-\dt}\,\,\text{ in $B_{2r}$. }
\end{equation}
If $N$ is the smallest natural number satisfying (8.2), then we easily obtain that
\begin{equation}N>\f{1}{\ln\,((1+\ep_0)^n\ap)}\,\ln\,\biggl(\f{|U^0_\tau|}{|B_{2r}|}\biggr).
\end{equation}
Thus by (8,2) and (8.3) we have that
$$\inf_{B_r} u\ge\inf_{B_{2r}}u\ge\tau\,\biggl(\f{|U^0_\tau|}{|B_{2r}|}\biggr)^{1/\sm}-\f{2b}{1-\dt}\,\,\text{ for $\sm=\f{\ln\,((1+\ep_0)^n\ap)}{\ln\dt}$, }$$
where $\dt$ and $\sm$ depend only on $n,s,t,p,q,\ld,\Ld,\Om$ and $\fa$.
This leads us to get that
\begin{equation}\f{|B_{2r}\cap\{u\ge\tau-2b/(1-\dt)\}|}{|B_{2r}|}=\f{|U^0_\tau|}{|B_{2r}|}\le\tau^{-\sm}\,\biggl(\,\inf_{B_r}u+\f{2b}{1-\dt}\biggr)^{\sm}.
\end{equation}
From stanard analysis, we see that
\begin{equation}\intavg_{B_{2r}}u^h\,dx=h\int_0^\iy\tau^{h-1}\,\f{|B_{2r}\cap\{u\ge\tau\}|}{|B_{2r}|}\,d\tau
\end{equation}
for any $h>0$.
Using (8.4) and (8.5), we have that
\begin{equation*}\begin{split}
\intavg_{B_{2r}}u^h\,dx&\le h\int_0^\iy\tau^{h-1}\,\f{|B_{2r}\cap\{u\ge\tau-2b/(1-\dt)\}|}{|B_{2r}|}\,d\tau \\
&\le h\int_0^{\tau_0}\tau^{h-1}\,d\tau+h\,\biggl(\,\inf_{B_r}u+\f{2b}{1-\dt}\biggr)^\sm\int_{\tau_0}^\iy\tau^{h-1-\sm}\,d\tau.
\end{split}\end{equation*}
Selecting $\tau_0=\inf_{B_r}u+2b/(1-\dt)$ and $h_0=\sm/2\in(0,1)$ in the above iequality, we obtain that
\begin{equation*}\intavg_{B_{2r}}u^{h_0}\,dx\le2\,\biggl(\,\inf_{B_r}u+\f{2b}{1-\dt}\biggr)^{h_0}
\le 2\,\biggl(\,\inf_{B_r}u+\f{8b}{3}\biggr)^{h_0}.
\end{equation*}
Thus it follows from H\"older's inequality with index $h_0/h\ge 1$ that
$$\biggl(\,\intavg_{B_{2r}}u^h\,dx\biggr)^{\f{1}{h}}\le \biggl(\,\intavg_{B_{2r}}u^{h_0}\,dx\biggr)^{\f{1}{h_0}}\le 2^{\f{1}{h_0}}\biggl(\,\inf_{B_r}u+\f{8b}{3}\biggr)$$ for any $h\in(0,h_0]$.
Therefore we are done. \qed

\section{ Nonlocal Harnack inequalities }

In this final section, we establish nonlocal Harnack inequalities for weak solutions of the nonlocal double phase equation (1.5) which is locally nonnegative in $\Om$.

\,

{\bf [Proof of Theorem 1.1.]} Take any $r>0$ satisfying $B_{2r}^0\subset B^0_R\subset\Om$. If $ps<n$ and $1<p\le q\le p_s^*$, then it follows from Theorem 5.3 and Lemma 1.3 that there is a constant $$c_0=c_0(n,s,t,p,q,\ld,\Ld,\Om,\fa)>0$$ such that
\begin{equation}\cT_r(u^+;x_0)\le c_0\,\biggl(\,\biggl[\sup_{B^0_r}u\biggr]^{p-1}+\biggl[\sup_{B^0_r}u\biggr]^{q-1} 
+\biggl(\f{r}{R}\biggr)^{ps}\cT_R(u^-;x_0)\,\biggr).
\end{equation} and 
\begin{equation}
\sup_{B^0_{r/2}}u\le\,\dt\,[\cT_r(u^+;x_0)]^{\f{1}{p-1}}+c_0\,\dt^{-\f{(p-1)n}{s p^2}}\biggl(\,\intavg_{B^0_{2r}}\Theta(u)\,dx\biggr)^{\f{1}{p}}.
\end{equation}
Let $\theta_0(u)$ be the value given by
$$\theta_0(u):=\biggl(\,\intavg_{B^0_{2r}}\Theta(u)\,dx\biggr)^{\f{1}{p}}.$$
We now consider two possible ocassions; (a) $\sup_{B^0_r}u\le1$ and (b) $\sup_{B^0_r}u>1$.

\,

[Case (a) $\sup_{B^0_r}u\le1$] In this case, by (9.1) we have that
\begin{equation}\cT_r(u^+;x_0)\le 2 c_0\,\biggl[\sup_{B^0_r}u\biggr]^{p-1}
+c_0\,\biggl(\f{r}{R}\biggr)^{ps}\cT_R(u^-;x_0).
\end{equation}
Thus it follows from (9.2) and (9.3) that
\begin{equation}\begin{split}
\sup_{B^0_{r/2}}u&\le c_0\,\dt^{-\f{(p-1)n}{s p^2}}\biggl(\,\intavg_{B^0_{2r}}\Theta(u)\,dx\biggr)^{\f{1}{p}} \\
&\qquad\qquad+\dt \,e_p(2c_0)^{\f{1}{p-1}}\,\sup_{B^0_r}u+e_p\, c_0^{\f{1}{p-1}}\dt\,\biggl(\f{r}{R}\biggr)^{\f{ps}{p-1}}
\cT^{\f{1}{p-1}}_R(u^-;x_0)
\end{split}\end{equation} for any $\dt\in(0,1/4)$,
where $e_p=\mathbbm{1}_{[2,\iy)}(p)+2^{\f{1}{p-1}}\mathbbm{1}_{(1,2)}(p)$ for $p>1$.

If $\theta_0(u)\le 1$, then we have that
\begin{equation}\begin{split}
\biggl(\,\intavg_{B^0_{2r}}\Theta(u)\,dx\biggr)^{\f{1}{p}}&\le \biggl(\,\intavg_{B^0_{2r}}u^p\,dx\biggr)^{\f{1}{p}}+\biggl(\,\intavg_{B^0_{2r}}\|\fa\|_\iy u^q\,dx\biggr)^{\f{1}{p}} \\
&\ls \biggl(\,\intavg_{B^0_{2r}}u^p\,dx\biggr)^{\f{1}{p}}+\biggl(\,\intavg_{B^0_{2r}} u^q\,dx\biggr)^{\f{1}{q}}. 
\end{split}\end{equation}
Thus it follows from (9.4) and (9.5) that
\begin{equation}\begin{split}
\sup_{B^0_{r/2}}u&\le c\,\psi(\dt)\,\biggl[\biggl(\,\intavg_{B^0_{2r}}u^p\,dx\biggr)^{\f{1}{p}}+\biggl(\,\intavg_{B^0_{2r}} u^q\,dx\biggr)^{\f{1}{q}}\biggr] \\
&\qquad+\dt \,e_p(2c_0)^{\f{1}{p-1}}\,\sup_{B^0_r}u+e_p\, c_0^{\f{1}{p-1}}\dt\,\biggl(\f{r}{R}\biggr)^{\f{ps}{p-1}}
\cT^{\f{1}{p-1}}_R(u^-;x_0)
\end{split}\end{equation}
where $\psi(\dt)$ is given by
$$\psi(\dt):=\dt^{-\f{(p-1)n}{s p^2}}.$$
By applying a covering argument with $1/2\le a<b\le 2$, we obtain that
\begin{equation}\begin{split}
\sup_{B^0_{ar}}u
&\le \f{c_*\,\psi(\dt)}{(b-a)^{\f{n}{p}}}\,\biggl[\biggl(\,\intavg_{B^0_{br}}u^p\,dx\biggr)^{\f{1}{p}}+\biggl(\,\intavg_{B^0_{br}} u^q\,dx\biggr)^{\f{1}{q}}\biggr] \\
&\qquad\qquad+c_*\dt\,\sup_{B^0_{br}}u+c_*\dt\,\biggl(\f{r}{R}\biggr)^{\f{ps}{p-1}}
\cT^{\f{1}{p-1}}_R(u^-;x_0) \\
&\le\f{c_*\,\psi(\dt)}{(b-a)^{\f{n}{p}}}\,\biggl[\biggl(\,\sup_{B^0_{br}}u\biggr)^{\f{p-h}{p}}\biggl(\,\intavg_{B^0_{br}}u^h\,dx\biggr)^{\f{1}{p}}+\biggl(\,\sup_{B^0_{br}}u\biggr)^{\f{q-h}{q}}\biggl(\,\intavg_{B^0_{br}} u^h\,dx\biggr)^{\f{1}{q}}\biggr] \\
&\qquad\qquad+c_*\dt\,\sup_{B^0_{br}}u+c_*\dt\,\biggl(\f{r}{R}\biggr)^{\f{ps}{p-1}}
\cT^{\f{1}{p-1}}_R(u^-;x_0)
\end{split}\end{equation}
Without loss of generality, we may assume that the constant $c_*$ in the above inequality satisfies $c_*>1$.
Taking $\dt=1/(4 c_*)$ in (9.7) and applying Young's inequality twice with indices $\f{p}{p-h},\f{p}{h},\vep=1/4$ and $\f{q}{q-h},\f{q}{h},\vep=1/4$ for $h\in(0,h_0]$ with the constant $h_0\in(0,1)$ in Theorem 1.6 makes it possible to obtain that
\begin{equation}\begin{split}
\sup_{B^0_{ar}}u&\le\f{1}{2}\sup_{B^0_{br}}u+\f{c_{h,p}+c_{h,q}}{(b-a)^{\f{n}{h}}}\,\biggl(\,\intavg_{B^0_{br}} u^h\,dx\biggr)^{\f{1}{h}}+\f{1}{4}\biggl(\f{r}{R}\biggr)^{\f{ps}{p-1}}
\cT^{\f{1}{p-1}}_R(u^-;x_0)
\end{split}\end{equation}
where $$c_{h,p}=\f{h}{p} 4^{\f{sp^2+(p-1)n}{sph}-1}\bigl(\f{p-h}{h}\bigr)^{\f{p-h}{h}}c_*^{\f{sp^2+(p-1)n}{sph}}$$ and
$$c_{h,q}=\f{h}{q} 4^{\f{tq^2+(q-1)n}{tqh}-1}\bigl(\f{q-h}{h}\bigr)^{\f{q-h}{h}}c_*^{\f{tq^2+(q-1)n}{tqh}}.$$
Employing Lemma 3.4 to (9.8), we have that
\begin{equation}\sup_{B^0_{\rho r}}u\le c\,\biggl[\f{c_{h,p}+c_{h,q}}{(2-\rho)^{\f{n}{h}}}\,\biggl(\,\intavg_{B^0_{2r}} u^h\,dx\biggr)^{\f{1}{h}}+\f{1}{4}\biggl(\f{r}{R}\biggr)^{\f{ps}{p-1}}
\cT^{\f{1}{p-1}}_R(u^-;x_0)\biggr]
\end{equation} for any $h\in(0,h_0]$ and any $\rho\in[1/2,2]$. Seleting $\rho=1$ in (9.9) and $h=h_0$ obtained in Theorem 1.6, we conclude that
\begin{equation}\sup_{B^0_r}u\ls\inf_{B^0_r}u+\biggl(\f{r}{R}\biggr)^{\f{ps}{p-1}}
\cT^{\f{1}{p-1}}_R(u^-;x_0)
\end{equation}
for any $r\in(0,R/2)$.

If $\theta_0(u)>1$, then we set $v=u/\theta_0(u)$. Then we claim that $\theta_0(v)\le 1$; indeed, we have the estimate
\begin{equation}\begin{split}
\intavg_{B^0_{2r}}\biggl(\,\f{u^p}{[\theta_0(u)]^p}+\|\fa\|_\iy\,\f{u^q}{[\theta_0(u)]^q}\biggr)\,dx
&\le\f{1}{[\theta_0(u)]^p}\intavg_{B^0_{2r}}(u^p+\|\fa\|_\iy u^q)\,dx=1.
\end{split}\end{equation}
Since $[\theta_0(u)]^{p-1}\le[\theta_0(u)]^{q-1}$, by (9.10) and (9.11) we have that
\begin{equation}\begin{split}
\f{1}{\theta_0(u)}\sup_{B^0_r}u&=\sup_{B^0_r}v\ls\inf_{B^0_r}v+\biggl(\f{r}{R}\biggr)^{\f{ps}{p-1}}
\cT^{\f{1}{p-1}}_R(v^-;x_0) \\
&\ls\f{1}{\theta_0(u)}\inf_{B^0_r}u+\f{1}{\theta_0(u)}\biggl(\f{r}{R}\biggr)^{\f{ps}{p-1}}
\cT^{\f{1}{p-1}}_R(u^-;x_0)
\end{split}\end{equation}
for any $r\in(0,R/2)$. Hence by (9.10) and (9.12) we conclude that
\begin{equation}\sup_{B^0_r}u\ls\inf_{B^0_r}u+\biggl(\f{r}{R}\biggr)^{\f{ps}{p-1}}
\cT^{\f{1}{p-1}}_R(u^-;x_0)
\end{equation}
for any $r\in(0,R/2)$.

[Case (b) $\sup_{B^0_r}u>1$] Set $M=\sup_{B^0_r}u$ and consider $v=u/M$. Then we see that 
\begin{equation*}\sup_{B^0_r}v\le 1.\end{equation*}
Since $M^{p-1}\le M^{q-1}$, applying $v$ to (9.13) yields that
\begin{equation*}\begin{split}
\f{1}{M}\sup_{B^0_r}u&=\sup_{B^0_r}v\ls\inf_{B^0_r}v+\biggl(\f{r}{R}\biggr)^{\f{ps}{p-1}}
\cT^{\f{1}{p-1}}_R(v^-;x_0) \\
&\ls\f{1}{M}\inf_{B^0_r}u+\f{1}{M}\biggl(\f{r}{R}\biggr)^{\f{ps}{p-1}}
\cT^{\f{1}{p-1}}_R(u^-;x_0) 
\end{split}\end{equation*}
for any $r\in(0,R/2)$. Also, even in this case, we conclude that
\begin{equation}\sup_{B^0_r}u\ls\inf_{B^0_r}u+\biggl(\f{r}{R}\biggr)^{\f{ps}{p-1}}
\cT^{\f{1}{p-1}}_R(u^-;x_0)
\end{equation}
for any $r\in(0,R/2)$. 

Finally, in case of $ps=n$ and $1<p\le q<\iy$, we can similarly obtain the inequality like (9.14) by using $h(\dt)=\dt^{-\f{q(p-1)}{p(q-p)}}$ given in Theorem 1.3. Hence we complete the proof. \qed

\,

\,

\noindent{\bf Data availability.} Data sharing not applicable to this article as no datasets were generated or used during this current study.

\,

\noindent{\bf Statements and Declarations(Competing Interests).} The author declares he has no competing interest.


\,

\end{document}